\documentclass{autart}
\usepackage{graphicx}      
\usepackage{natbib}        
\usepackage{epsfig} 
\usepackage{amsmath} 
\usepackage{amssymb}  
\usepackage{epstopdf}
\usepackage{color}
\usepackage{xcolor}
\usepackage{arydshln}
\usepackage{subfigure}
\usepackage{cases}
\usepackage{mathptmx} 
\usepackage{times}
\usepackage{subeqnarray}
\usepackage{appendix}

\newcommand{\dn}{\mathbf{d}}
\newcommand{\tdn}{\tilde{\mathbf{d}}}
\newcommand{\R}{\mathbb{R}}
\newcommand{\diag}{\mathrm{diag}}

 \newtheorem{Def}{Definition}
  \newtheorem{remark}{Remark}
 \newtheorem{Pro}{Proposition}
\begin{document}
\graphicspath{{figure/}}
\begin{frontmatter}

\vspace{-1cm}
\title{Homogeneous Distributed Observers for Quasilinear Systems}
 \vspace{-1cm}

\thanks[footnoteinfo]{Corresponding author Siyuan Wang.}


\author[First]{Min Li}\ead{min.li@inria.fr},
\author[First]{Andrey Polyakov}\ead{andrey.polyakov@inria.fr},
\author[Second]{Siyuan Wang}\ead{siyuanwang0603@buaa.edu.cn},
\author[First]{Gang Zheng}\ead{gang.zheng@inria.fr}
\vspace{-0.2cm}
\address[First]{Inria, Centrale Lille, Universit\'{e} de Lille, CNRS, CRISTAL, France.}\vspace{-0.2cm}

\address[Second]{State Key Laboratory of Virtual Reality Technology and Systems, Beihang University, Beijing, China.}\vspace{-0.8cm}




\begin{keyword}
\!\!\!Distributed observer, H{\"o}lder condition, generalized homogeneity, finite-time stability, fixed-time stability, input-to-state stability.\vspace{-0.2cm}
\end{keyword}

\vspace{-0.1cm}
\begin{abstract}
\vspace{-0.8cm}

The problem of finite/fixed-time cooperative state estimation is considered for a class of quasilinear systems with nonlinearities satisfying a H{\"o}lder condition. A strongly connected nonlinear distributed observer is designed under the assumption of global observability.  By proper parameter tuning with linear matrix inequalities, the observer error equation possesses finite/fixed-time stability in the perturbation-free case and input-to-state stability with respect to bounded perturbations. Numerical simulations are performed to validate this design. \vspace{-0.6cm}
\end{abstract}

 \end{frontmatter}

\vspace{-1cm}
\section{Introduction}\label{Intro}
\vspace{-0.3cm}

Traditional Luenberger state observers, which rely on global outputs, cannot be effectively designed for some Cyber-Physical Systems with geographically dispersed sensors/measurements (see \textit{e.g.,} \cite{mitra2018distributed}). Facing this challenge, a cooperative state estimation technique that incorporates multiple observers to design a distributed observer for the linear plant is developed (\cite{carli2008distributed,park2012necessary,mitra2018distributed,kim2016distributed,Liu2017CooperativeStabilization,han2018simple}).
\cite{BATTILOTTI2019562,xu2021distributed,wu2021design} focus on distributed observer designs for nonlinear plants with global Lipschitz nonlinearities.  To the best of the authors' knowledge, distributed observers for nonlinear systems, which do not satisfy global Lipschitz conditions, are not developed yet. However, centralized observers for plants modeled by nonlinear systems satisfying  H$\Ddot{\mathrm{o}}$lder condition can be found in the recent literature (see, \textit{e.g.,}  \cite{du2013recursive, bernard2017observers,bernard2022observer}).
\vspace{-0.2cm}

\cite{BATTILOTTI2019562,xu2021distributed,wu2021design} primarily design the asymptotic distributed observer such that the state estimation is guaranteed as time goes to infinity. In such cases, ensuring the performance/quality of state estimation is challenging. To deal with this issue, some studies are devoted to the tuning of the convergence rate of distributed observers. 
\cite{han2018simple} propose state estimates with guaranteed exponential convergence rates. 
\vspace{-0.2cm}


The finite-time stability is a popular approach addressing the non-asymptotic control/observation problem. Being known since 1960s (\cite{Fuller1960:IFAC,Korobov1979:DAN,Haimo1986:SIAM,bhat2000finite}), and has received extensive attention in the last two decades (\cite{EfimovPolyakov2021:Survey}). A system is said to be finite-time stable if it is Lyapunov stable and its states reach zeros at a finite instant of time (which may depend on the initial state of the system). To achieve finite-time stability, there are two main approaches, one is based on the finite-time Lyapunov function theory as \cite{bhat2000finite}; the other one is by utilizing the property of the homogeneous system, $\it{i.e.,}$ the asymptotic stability of such a system with a negative homogeneity degree implies its finite-time stability (\cite{bhat2005geometric}). The fixed-time stability  (\cite{polyakov2011nonlinear}),
is a further development of the finite-time stability concept, which assumes the uniform boundedness of the settling-time function
for all admissible initial conditions. 
\vspace{-0.2cm}

The integration of the finite/fixed-time techniques to centralized observer design has been actively researched, see $\textit{e.g.,}$ \cite{engel2002continuous,andrieu2008homogeneous,perruquetti2008finite,lopez2018finite,kitsos2021high}. However, for the distributed observer, few results are reported. In particular, \cite{ortega2020distributed,silm2018note, silm2020simple} focus on the finite-time distributed observer, and \cite{Ge2023TAC} present the sole result on the fixed-time distributed observer, utilizing a kernel-based design approach. It's worth noting that all these studies concentrate on LTI plants and exhibit a notable gap in addressing the cooperative estimation of nonlinear plants. 
Based on the review of existing results, it is essential to develop a distributed observer admitting finite/fixed-time cooperative estimation for a class of nonlinear plants yield H$\Ddot{\mathrm{o}}$lder conditions. 
\vspace{-0.2cm}

Building on the preceding discussion, this paper assumes the observed plant follows a quasilinear model, a fundamental nonlinear model containing a linear component and a nonlinear component yields H$\Ddot{\mathrm{o}}$lder conditions. The distributed observer is primarily designed for the linear component by upgrading the classical linear distributed observers to a homogeneous one, which allows finite-time stability of estimation errors with a practical parameter tuning procedure similar to the linear case.  
Based on the homogeneous design, a slight modification allows the fixed-time cooperative estimation. In both finite/fixed-time designs, the robustness of the estimation with respect to bounded perturbations such as the noise of states and measurement outputs is guaranteed. Finally,  the observer design is completed by characterizing the H$\Ddot{\mathrm{o}}$lder condition that the observed plant follows such that all aforementioned properties are preserved. 
\vspace{-0.2cm}

This paper is organized as follows. Section \ref{Pre} gives some basic knowledge of graph theory, stability definitions, and generalized homogeneity. The problem to be studied is formulated in Section \ref{PS}. The basic idea of the observer design is proposed in Section \ref{BI}. The main results are presented in Section \ref{MR}. Finally, in Section \ref{SR}, the effectiveness of the proposed observer is illustrated by some numerical simulations.
\vspace{-0.6cm}

\textit{Notation}.
$\R$ is the field of reals and  $\mathbb{R}_{\geq0}\!\!=\!\!\{x\!\!\in\!\! \R\!\!:\!\! x\!\!\geq\!\! 0\}$;
$\mathbb{N}_+$ is the set of positive integers;
given $n\!\!\in\!\! \mathbb{N}_+$ a sequence of positive integers $1,\dots,n$ is denoted as $\overline{1,n}$; 
$\R^n$ and $\R^{n\times n}$ denote the  $n\!\!\times\!\! 1$ real vector and the $n\!\!\times\!\! n$ real matrix, respectively;  $I_n$ is the $n\!\!\times\!\! n$ identity matrix; $\mathbf{1}_n\!\!\in\!\! \R^n$ is the vector whose components are all ones; $\mathbf{0}$ is the zero elements (\textit{e.g.,} zero matrices, zero vectors, zero functions, etc); $\diag\{\sigma_i\}_{i=1}^n$ is the (block) diagonal matrix with elements $\sigma_i$ of a proper dimension;  $\{\sigma_{ij}\}$ is a $n_1\!\!\times \!\!n_2$ (block) matrix of a proper dimension whose elements are $\sigma_{ij}$, $i\!\!=\!\!\overline{1,n_1}$, $j\!\!=\!\!\overline{1,n_2}$, $n_1,n_2\!\!\in\!\!\mathbb{N}_+$; 
$P\!\!\succ\!\! 0(\textit{resp.,} \prec \!\!0)$ for $P\!\!\in\! \!\R^{n\times n}$ means that the matrix $P$ is symmetric and positive (\textit{resp.,} negative) definite;
$\|\cdot\|$ is a norm in $\R^n$; $|\cdot|$ is the Euclidean norm in $\R^n$; $\|x\|_{\infty}\!\!=\!\!\max_{i=\overline{1,n}}\!|x_i|$ for $x\!\!=\!\!(x_1,...,x_n)^{\top}\!\!\!\in\!\! \R^n$;
$\exp(Q)\!\!=\!\!\sum_{i=0}^{\infty} \frac{Q^i}{i!}$  for $Q\!\!\in\!\! \R^{n\times n}$ with $n\!\!\in\!\! \mathbb{N}_+$;
$\otimes $ denotes the Kronecker product;
$C^p(X,Y)$, $p\!\!\in\!\!\mathbb{N}_+\!\!\cup\!\!\{0\}$ is the class of functions $X\!\!\to\!\!Y$ which are continuously differentiable at least up to the $p$-order, where $X$ and $Y$ are the subsets of normed vector space;
a function $\alpha\!\!\in\!\!C(\mathbb{R}_{\geq0},\mathbb{R}_{\geq0})$
is of the class $\mathcal{K}$ if it is strictly increasing with $\alpha(0)\!\!=\!\!0$, 
a function $\beta\!\!\in\!\! C(\mathbb{R}_{\geq0}\!\!\times\!\mathbb{R}_{\geq0},\mathbb{R}_{\geq0})$ 
is of the class $\mathcal{KL}$ if $\beta(\cdot,t)\!\!\in\!\!\mathcal{K}$ for each fixed $t$, and, 
for each fixed $\varrho$, the function $t\mapsto \beta(\varrho,t)$ is strictly decreasing  to zero;
$L^\infty(\R_{\geq0},\R^d)$ is a Lebesgue space of measurable uniformly essentially bounded functions 
$q\!\!:\!\!\R_{\geq 0}\!\!\to\!\!\R^d$ with the norm $\|q\|_{L^\infty}\!\!:=\!\!\mathrm{ess}\sup_{t\geq0}\!\|q(t)\|_\infty$;
${B}(r)\!\!=\!\!\{x\!\!\in\!\! \R^n\!\!:\!\! \|x\|\!\!<\!\! r\}$ denotes the open ball of
radius $r\!\!>\!\!0$ centered at the origin.

\vspace{-0.4cm}
\section{Preliminaries}\label{Pre}
\vspace{-0.3cm}
\subsection{Graph Theory}
\vspace{-0.3cm}
A fixed directed graph $\mathcal{G}$ is usually characterized by a $\it{node\;set}\;\mathcal{V}$, an $\it{edge\;set}\;\mathcal{E}$ and an $\it{adjacency\;matrix}\;\mathcal{A}$. To be specific, let a graph with $N\!\!\in\!\! \mathbb{N}_+$ nodes, the node set $\mathcal{V}\!\!=\!\!\{1,\dots,N\}$ contains all the nodes at the graph labeled by $i\!\!=\!\!\overline{1,N}$; the edge set $\mathcal{E}\!\!=\!\!\{(i,j)|i,j\!\!\in\!\! \mathcal{V}\}$, we have $(i,j)\!\!\in\!\! \mathcal{E}$ if the node $j$ is able to transfer its local information to the node $i$, the number of incoming edges of node $i$ is denoted by $n_i$; $\mathcal{A}\!\!=\!\!\{a_{ij}\}\!\!\in\!\!\R^{N\times N },\;a_{ij}\!\!\in\!\!\R$, where $a_{ij}\!\!=\!\!1$ if $(i,j)\!\! \in\!\! \mathcal{E}$ and $a_{ij}\!\!=\!\!0$ otherwise.
A $\textit{directed path}$ from the node $i$ to the node $j$ is a sequence of nodes $i_0,\ldots, i_s$, where $i_0\!\!=\!\!i$, $i_s\!\!=\!\!j$ and $(i_{\kappa+1},i_\kappa)\!\!\in\!\!\mathcal{E},\;\kappa\!\!=\!\!\overline{1,s\!-\!1}$.
The graph $\mathcal{G}$ is $\textit{strongly connected}$ if there exists at least one directed path between each pair of the nodes. In this paper, we assume self-loops are excluded, \textit{i.e.,} $(i,i) \!\!\not\in\!\! \mathcal{E}$.
The $\textit{Laplacian matrix}$ associated to graph $\mathcal{G}$ is defined as $\mathcal{L}\!\!=\!\!\{l_{ij}\}\!\!\in\!\!\R^{N\times N },\;l_{ij}\!\!\in\!\!\R$, where $l_{ij}\!\!=\!\!-a_{ij}$ if $i\!\!\neq\!\! j$ and $l_{ij}\!\!=\!\!\sum\nolimits^{N}_{k=1}\!a_{ik}$ if $i\!\!=\!\! j$.
For the strongly connected graph, the associated Laplacian matrix has a simple $0$ eigenvalue while all others have positive real parts, and the associated left $0$-eigenvector $\zeta^\top\!\!\!\!=\!\!(\zeta_1,\dots,\zeta_N)$  yields $\zeta_i\!\!>\!\!0$, $i\!\!=\!\!\overline{1,N}$.
\vspace{-0.3cm}\begin{lem}\label{lem_laplacian}[\cite{lewis2013cooperative}]
{\it{Let $\mathcal{L}\!\!\in\!\!\R^{N\times N}$ be the Laplacian matrix corresponding to a strongly connected graph $\mathcal{G}$, then $\mathcal{L}$ can be similarly transformed by $T\!\!\in\!\!\R^{N\times N}$ and its inverse, which is given as
\vspace{-0.3cm}
\begin{equation*}
\begin{aligned}
  T^{-1}\!\!{\mathcal{L}}T\!\!=\!\!\left(\begin{smallmatrix}
    0 & \textbf{0}^{\top}_{N\!-\!1}\\
    \textbf{0}_{N\!-\!1} & \Delta
\end{smallmatrix}\right),
\end{aligned}
\vspace{-0.3cm}\end{equation*}
where $\Delta$ is a matrix whose eigenvalues correspond to the nonzero eigenvalues of ${\mathcal{L}}$, $T^{-1}\!\!=\!\!(\zeta,Y_2)^\top$, $T\!\!=\!\!(\mathbf{1}_N, Y_1)$, with $\zeta$ being the left 0-eigenvector of $\mathcal{L}$ yields $\zeta^\top\!\textbf{1}_N\!\!=\!\!1$ by normalization, $Y_1$ (\textit{resp.,} $Y_2$) is composed of the right (\textit{resp.,} left)-eigenvectors of ${\mathcal{L}}$ associated to the nonzero eigenvalues, yields $Y_2^{\top}Y_1\!\!=\!\!I_{N-1}$. }}
\end{lem}

\vspace{-0.3cm}
\subsection{Generalized Homogeneity}
\vspace{-0.3cm}

Homogeneity is an invariance of an object with respect to a class of transformations called dilations (\cite{Kawski1991:ADCS}, \cite{Zubov1957:IVM}). Choosing a proper dilation group $\mathbf{d}(s),\;s\!\in\! \mathbb{R}$ is vital for the homogeneity-based analysis, $\mathbf{d}(s)$ is supposed to satisfy the limit property: $\lim_{s\rightarrow \pm \infty}\!\|\mathbf{d}(s)x\|\!\!=\!\!\exp(\pm \infty)$ for all $ x\!\!\in\!\!\mathbb{R}^n\!\setminus\!\!\{\mathbf{0}\}$. 
A dilation $\mathbf{d}(s)$ is \textit{monotone} with respect to the norm $\|\cdot\|$ if the function $s\!\!\mapsto\! \!\|\mathbf{d}(s)x\|$ is strictly
increasing for any $x\!\!\in\!\!\mathbb{R}^n\!\setminus\!\!\{\mathbf{0}\}$.
This work uses the $\it{linear\;dilation}$ (\cite{Polyakov2020:Book}) which is defined as
$\mathbf{d}(s)\!\!=\!\!\exp(s G_{\dn})$, $s\!\!\in\!\! \mathbb{R}$,
where $G_{\dn}\!\!\in\!\! \mathbb{R}^{n\times n}$ is an anti-Hurwitz matrix known as the generator of the dilation.
The case $G_{\dn}\!\!=\!\!I_n$ corresponds to the so-called standard (or Euler) dilation being a multiplication of a vector by a positive scalar $\exp(s)$. Any other dilation is called  \textit{generalized} and the corresponding homogeneity is known as generalized homogeneity (\cite{Zubov1957:IVM}, \cite{Khomenuk1961:IVM}, \cite{Polyakov2020:Book}). 
\vspace{-0.3cm}
%

\begin{Def}\label{def:homogeneity}[Homogeneous Vector Fields (\cite{Kawski1991:ADCS,Polyakov2020:Book})]
{\it{	A vector field $f\!\!:\!\! \mathbb{R}^n \!\!\to\!\! \mathbb{R}^n$ 
is $\dn$-homogeneous if there exists  $\mu \!\!\in\!\! \mathbb{R}$ such that 
	$f(\dn(s) x)\!\!=\!\!\exp(\mu s) \dn(s) f(x)$, $s\!\!\in\!\! \mathbb{R}$,  $x \!\!\in\!\! \mathbb{R}^n$, 
where $\dn$ is a dilation and $\mu$ is the so-called  homogeneity degree.}}
\end{Def}
\vspace{-0.3cm}
\begin{Def}\label{def:cano_norm}[Canonical Homogeneous Norm (\cite{Polyakov2020:Book})]
{\it{	The function $\|x\|_{\mathbf{d}}\!\!:\!\!\mathbb{R}^n\!\!\to\!\! \R_{\geq 0}$ defined as 
$\|\mathbf{0}\|_{\dn}\!\!=\!\!0,$  
\vspace{-0.3cm}\begin{equation}\label{def:hnorm}
		\|x\|_\mathbf{d}\!\!=\!\!{\rm exp}(s_x),\;\; s_x\!\in\! \mathbb{R}\!:\!\|\mathbf{d}(-s_x)x\|\!\!=\!\!1,\;\; x\!\in\! \mathbb{R}^n\!\!\setminus\!\!\{\mathbf{0}\}
  \vspace{-0.3cm}\end{equation}
  is called the canonical homogeneous norm in $\mathbb{R}^n$, with $\mathbf{d}$ being a linear monotone dilation.}}
\end{Def}\vspace{-0.3cm}
If the linear dilation $\dn(s)\!\!=\!\!\exp(sG_{\dn})$ is monotone with respect to the norm $\|x\|\!\!=\!\!\sqrt{x^{\top}\!\!Px}$  then (\cite{Polyakov2018:RNC})
\vspace{-0.3cm}\begin{equation}\label{eq:procannorm}
		\tfrac{\partial \|x\|_{\dn}}{\partial x}\!\!=\!\!\tfrac{\|x\|_{\dn}x^{\top}\!\dn^{\top}\!(-\ln \|x\|_{\dn}) P\dn(-\ln \|x\|_{\dn})}{x^{\top}\!\dn^{\top}\!(-\ln \|x\|_{\dn})PG_{\dn}\dn(-\ln \|x\|_{\dn})x}, \;\;\; x\!\!\in\!\!\R^n\backslash\{\mathbf{0}\},
	\vspace{-0.2cm}\end{equation}
 \vspace{-0.3cm}\begin{equation}\label{eq:hom_norm_est}
 \min\{\|x\|^\alpha_{\dn},\|x\|_{\dn}^{\beta}\}\!\!\leq\!\!\|x\|\!\!\leq\!\!\max\{\|x\|^\alpha_{\dn},\|x\|^\beta_{\dn}\},
 \vspace{-0.05cm}\end{equation}
 where $\alpha$ and $\beta$ are the maximal and minimal eigenvalue of the matrix $PG_{\dn}\!\!+\!\!G_{\dn}^{\top}\!\!P\!\!\succ\!\! 0$.

\vspace{-0.3cm}
\begin{lem}\label{lem:homo_fin}[\cite{nakamura2002smooth}]
{\it{Let $f\!\!:\!\!\R^n\!\!\to\!\!\R^n$ be a $\dn$-homogeneous vector
field of degree $\mu$. If the system $\dot{x}\!\!=\!\!f(x)$ is globally uniformly asymptotically stable then it is globally uniformly finite-time stable\footnote{The system $\dot{x}(t)\!\!=\!\!f(x(t))$, $t\!\!>\!\!0$, $x(0)\!\!=\!\!x_0\!\!\in\!\! \R^n,$ is globally uniformly finite-time stable (\cite{Orlov2005:SIAM_JCO}) if it is Lyapunov stable and there exists a locally bounded settling-time function $T(x_0)$, $T\!\!:\!\!\mathbb{R}^{n}\!\!\to\!\!\mathbb{R}$ such that $\|x(t)\|\!\!=\!\!0$ for all $t\!\!\geq\!\! T(x_0)$. } if $\mu\!\!<\!\!0$ and globally uniformly nearly fixed-time stable\footnote{The system $\dot{x}(t)\!\!=\!\!f(x(t))$, $t\!\!>\!\!0$, $x(0)\!\!=\!\!x_0\!\!\in\!\! \R^n,$ is globally uniformly nearly fixed-time stable (\cite{EfimovPolyakov2021:Survey}) if it is Lyapunov stable and for all $r\!\!>\!\!0$, there exists $T_r\!\!>\!\!0\!:\!\|x(t)\|\!\!<\!\!r$ for all $t\!\!\geq\!\! T_r$, where $T_r$ is independent of $x_0$. } if $\mu\!\!>\!\!0$.}}
\end{lem}
\vspace{-0.2cm}
\section{Problem Statement}\label{PS}
\vspace{-0.3cm}
Consider a plant described by a quasilinear system:
\vspace{-0.3cm}
\begin{equation}\label{eq:LTI_plant}
\dot x(t) \!\!=\!\! Ax(t)\!+\!Bu(t)\!+\!{\gamma(t,x)}\!+\!q_x(t),\;\;\; t\!\!>\!\!0,
\vspace{-0.3cm}
\end{equation}
where $x\!\!\in\!\!\mathbb{R}^n$ is the plant state, $A\!\!\in\!\! \mathbb{R}^{n\!\times \!n}$, $B\!\!\in\!\!\mathbb{R}^{n\!\times\! m}$, $u\!\!\in\!\! \mathbb{R}^{m}$ is the control input,  $\gamma\!\!\in\!\!C(\R\!\!\times\! \!\R^n,\R^n$) is a nonlinear function satisfying a H\"older condition to be given below, $q_x\!\!\in\!\!L^\infty(\mathbb{R},\mathbb{R}^{n})$ is unknown additive perturbation.




\vspace{-0.3cm}
A set of distributed sensors generates local measurements (outputs) of the plant \eqref{eq:LTI_plant}:
\vspace{-0.3cm}
\begin{equation}\label{eq:sensor_output}
y_i(t)\!\!=\!\! C_ix(t)\!+\!q_{y,i}(t),\;\;i\!=\!\overline{1,N}
\vspace{-0.3cm}
\end{equation}
where $y_i\!\!\in\!\!\mathbb{R}^{p_i}$, $C_i\!\!\in\!\!\mathbb{R}^{p_i\times n}$, $q_{y,i}\!\!\in\!\!L^\infty(\mathbb{R},\mathbb{R}^{p_i})$ models a measurement noise. A  topology of the sensor network is defined by a fixed directed graph $\mathcal{G}\!\!=\!\!\{\mathcal{V},\mathcal{E},\mathcal{A}\}$. Let the matrix $C\!\!=\!\!(C^{\top}_1\!\!,\dots,C^{\top}_N)^{\top}\!\!\in\! \mathbb{R}^{p\times n}$, $p\!\!=\!\!\sum^{N}_{i=1}\!p_i$ be a collection of $C_i$. 


\vspace{-0.3cm}
\begin{assum}\label{assum2}
 {\it{The pair $(A,C)$ is observable.}}
\vspace{-0.3cm}\end{assum}
This assumption is necessary for the existence of any finite/fixed-time observer.
Recall \cite{Polyakov2020:Book,zimenko2020robust} that for observable pair $(A,C)$  the following algebraic equations  
\vspace{-0.3cm}
\begin{equation}\label{eq:G_0}
	G_0A\!-\! AG_0\!+\!Y_0C\!\!=\!\!A, \;\;  CG_0\!\!=\!\!\mathbf{0},
	\vspace{-0.3cm}
 \end{equation}
 always have a solution $(Y_0,G_0)\!\!\in\!\!(\R^{p\times n},\R^{n\times n})$ 
such that $I_n\!+\!G_0$ is invertible  and the matrix $A_0\!\!=\!\!A\!+\!L_0C$  is nilpotent for 
$L_0\!\!=\!\!(I_n\!+\!G_0)^{-1}Y_0$. The latter implies the linear vector field $x\!\mapsto\!\!A_0x$, $x\!\!\in\!\!\R^n$ to be homogeneous. 
\vspace{-0.2cm}

Below we distinguish two classes of systems under consideration:
\textit{a system with $L_0\!\!=\!\!\textbf{0}$ and a system with $L_0\!\!\neq\!\! \textbf{0}$}. We consider the matrix $L_0$ as a parameter induced by the plant model, which impacts the observer design algorithm. 
\vspace{-0.2cm}

This work deals with the design of a  distributed observer composed of a set of observers with dynamics:
\vspace{-0.3cm}
\begin{equation}\label{eq:goalobserver}
  \dot {\hat x}_i\!\!=\!\!\xi_i(\hat x_i,y_i,\hat x_{j_1},\dots,\hat x_{j_{n_i}}),\;\;\;i\!\!= \!\!\overline{1,N}\vspace{-0.3cm}
\end{equation}
 which cooperatively reconstruct the state of (\ref{eq:LTI_plant}) in a finite/fixed time in the disturbance-free case, and which is robust (in ISS\footnote{Input-to-State Stable. The system $\dot x(t)\!\!=\!\!f(x(t),q(t))$, $x\!\!\in\!\! \R^n$, $q\!\!\in\!\! \R^d$, $t\!\!>\!\!0$, $x(0)\!\!=\!\!x_0$ is  ISS (\cite{Sontag1989:TAC}) if there exist $\beta\!\!\in\!\! \mathcal{KL}\!\!$ and $\alpha\!\!\in\!\!\mathcal{K}$ such that
 \vspace{-0.3cm}\begin{equation*}
	\|x(t)\|\!\!\leq\!\! \beta\left(\|x_0\|,t\right)\!\!+\!\!\alpha(\|q\|_{L^{\infty}[0,t)}),
	\vspace{-0.3cm}
 \end{equation*}
	 for any $x_0\!\!\in\!\! \R^n$ and any $q\!\!\in\!\! L^{\infty}(\R,\R^d)$.} sense) with respect to perturbations $q_x$ and $q_y$,
 where $\hat x_i\!\!\in\!\!\R^n$ is the state of the $i$-th observer, $j_k\!\!\in\!\! \mathcal{V}:(i,j_k)\!\!\in\!\! \mathcal{E}$, $\xi_i\!\!:\!\!\R^{n}\!\!\times\!\!\R^{p_i}\!\!\times\!\!\R^{n\cdot n_i}\!\!\to\!\!\R^n$
 and $q_{y}\!\!=\!\!(q_{y,1}^\top,\dots,q_{y,N}^\top)^\top\!\!\!\!\in\!\!L^\infty(\mathbb{R},\mathbb{R}^{p})$.
\vspace{-0.2cm}
Our \textit{first aim} is to design the distributed observer \eqref{eq:goalobserver} such that the estimation error equation
\vspace{-0.3cm}
\begin{equation}\label{eq_aug}\begin{aligned}
	\dot e\!\!=&f(e,q), \;\; f\!:\!\mathbb{R}^{ Nn}\!\!\times\!\!\mathbb{R}^{n+p}\!\!\to\!\!\mathbb{R}^{Nn},\\
e\!\!=&((\hat{x}_1\!-\!x)^{\top}\!\!, \dots, (\hat{x}_{N}\!-\!x)^{\top} )^{\top}\!\!,\;\; q\!\!=\!\!(q_x^{\top},q_y^{\top})^{\top},
\end{aligned}\vspace{-0.3cm}
\end{equation}
for $\gamma\!\!=\!\!\textbf{0}$ and $L_0\!\!=\!\!\textbf{0}$, has the following properties: 
\vspace{-0.2cm}

$\bullet$ there exists a linear dilation $\tilde\dn$ in $\R^{Nn}$ such that the vector field $f(\cdot,\textbf{0})$ is $\tilde\dn$-homogeneous of degree $\mu$;
\vspace{-0.2cm}

$\bullet$ the unperturbed ($q\!\!=\!\!\mathbf{0}$) error equation is globally uniformly finite-time (\textit{resp.,} exponentially or nearly fixed-time) 
stable for $\mu\!\!<\!\!0$ 
(\textit{resp.,} $\mu\!\!=\!\!0$ or $\mu\!\!>\!\!0$);
\vspace{-0.2cm}

$\bullet$ the error equation is ISS with respect to $q\!\!\in\!\!L^\infty(\mathbb{R},\mathbb{R}^{n+p})$.
\vspace{-0.2cm}

The \textit{second aim} is to modify the homogeneous observer design 
such that for $\gamma\!\!=\!\!\textbf{0}$ and 
$L_0\!\!\neq\!\!\textbf{0}$  the error equation \eqref{eq_aug} is
\vspace{-0.2cm}

$\bullet$ globally uniformly fixed-time stable\footnote{The system $\dot{x}(t)\!\!=\!\!f(x(t))$, $t\!\!>\!\!0$, $x(0)\!\!=\!\!x_0\!\in\! \R^n,$ is globally uniformly fixed-time stable (\cite{Polyakov2020:Book}) if it is globally uniformly finite-time stable and the settling-time function $T(x_0)$ is globally bounded. } provided $q\!\!=\!\!\textbf{0}$;
\vspace{-0.2cm}

$\bullet$ ISS with respect to $q\!\!\in\!\!L^\infty(\mathbb{R},\mathbb{R}^{n+p})$.
\vspace{-0.2cm}

Our \textit{final aim} is to design observers for $\gamma\!\!\neq\!\!\textbf{0}$ and to characterize a class of admissible nonlinearities allowing the error equation to preserve the properties mentioned above.

\vspace{-0.4cm}
\section{Basic Idea of the Observer Design}\label{BI}
\vspace{-0.3cm} 
For $\gamma\!\!=\!\! \textbf{0}$, inspired by \cite{kim2016distributed,Liu2017CooperativeStabilization,han2018simple}, let a linear distributed observer:
\vspace{-0.3cm}\begin{equation}\label{eq_liner_observer}\begin{aligned}
  &\dot{\hat{x}}_i \!\!=\!\!A\hat{x}_i\!+\!Bu\!+\!H_i\omega_i\!+\!\nu \theta_i,\;\;i\!\!=\!\!\overline{1,N}
\end{aligned}\vspace{-0.3cm}
\end{equation}
where $\omega_i\!\!=\!\!C_i\hat{x}_i\!-\!y_i\!\!=\!\!C_ie_i\!-\!q_{y,i}$, $\theta_i\!\!=\!\!\textstyle{\sum}^{N}_{j=1}a_{ij}(\hat{x}_j\!-\!\hat{x}_i)\!\!=\!\!-(\mathcal{L}_i\!\otimes\! I_n)e$ with $\mathcal{L}_i$ being the $i_{\mathrm{th}}$ row of the Laplacian matrix $\mathcal{L}$. 
The global estimation error equation is
\vspace{-0.3cm}\begin{equation}\label{eq:lin_error_dy}\begin{aligned}
  &\dot{e} \!\!=\!\!(\tilde A\!+\!\tilde{H}\tilde{C}\!-\!\nu\mathcal{L}\!\otimes\! I_n)e\!-\!\tilde{H}q_y\!-\!\textbf{1}_N\!\otimes \!q_x,\\&\tilde A\!\!=\!\!I_N\!\otimes\! A,\;\tilde H\!\!=\!\!\diag\{H_i\}_{i=1}^N,\;\tilde C\!\!=\!\!\diag\{C_i\}_{i=1}^N.
\end{aligned}
\vspace{-0.3cm}\end{equation}
If $(A,C)$ is observable, then the LMI
 \vspace{-0.3cm}\begin{equation}\label{LMI_Global_Observ}
P_a\!\!\succ\!\!0,\;\;\;P_aA\!+\!A^\top\!\!P_a\!+\!YC\!+\!C^{\top}\!\!Y^{\top}\!\!\!+\!2\rho P_a\!\!\prec\!\!0
		\vspace{-0.3cm}
  \end{equation}
is feasible for $\rho\!\!>\!\!0$ (see \cite{Boyd_etall1994:Book}) and $A\!+\!\bar H C$ is Hurwitz, where $\bar H\!\!=\!\! P^{-1}_aY\!\!\in\!\! \mathbb{R}^{n\times p}$,  $P_a\!\!\in\!\! \mathbb{R}^{n\times n}$, $Y\!\!\in\!\! \mathbb{R}^{n\times p}$.


\vspace{-0.3cm}
\begin{lem}\label{thm_linear}[refined from \cite{kim2016distributed,Liu2017CooperativeStabilization,han2018simple}]
{\itshape Let Assumption \ref{assum2} be fulfilled. Let the linear distributed observer \eqref{eq_liner_observer} operate under a strongly connected graph $\mathcal{G}$ with Laplacian matrix $\mathcal{L}$. Let $P_a\!\!\in\!\!\mathbb{R}^{n\times n}$ and $\bar H\!\!=\!\!P^{-1}_aY\!\!=\!\! (\bar H_1,\dots, \bar H_{N})\!\!\in\!\! \mathbb{R}^{n\times p}$, $\bar H_i\!\!\in\!\! \mathbb{R}^{n\times p_i}$ be defined by solving  \eqref{LMI_Global_Observ} and 
\vspace{-0.3cm}\begin{equation}\label{LMI:DeltaP}
   P_a\!\!\succ\!\!0,\;\;\;(\Delta^\top\!\!\!+\!\Delta)\!\otimes\!P_a\!\succ\!0, 
\vspace{-0.3cm}\end{equation}
where $\Delta$ is a matrix whose eigenvalues
correspond to the nonzero eigenvalues of $\mathcal{L}$.
 Let $H_i\!\!=\!\!\frac{\bar{H}_i}{\zeta_i}$, where $\zeta_i\!\!>\!\!0$ is the element of $\zeta$ defined as the left $0$-eigenvector  of $\mathcal{L}$ which yields $\zeta^\top\!\!\textbf{1}_N\!\!=\!\!1$. Then for a  large enough $\nu\!\!>\!\!0$ 
then the error equation \eqref{eq:lin_error_dy} is
globally uniformly exponentially stable  if  $q\!\!=\!\!\mathbf{0}$ and %
ISS with respect to $q\!\!\in\!\!L^\infty(\mathbb{R},\mathbb{R}^{n+p})$.
            }
\end{lem}
\vspace{-0.6cm}
\begin{pf}
Let $q\!\!=\!\!\textbf{0}$. Let $\eta\!\!=\!\!(T^{-1}\!\otimes\! I_n)e$, where $T\!\!=\!\!(\mathbf{1}_N,Y_1)$
and $T^{-1}\!\!=\!\!(\zeta,Y_2)^{\top}$ are given in Lemma \ref{lem_laplacian}. 
Since 
\vspace{-0.3cm}\begin{equation*}
\begin{aligned}
&(\zeta^{\top}\!\!\!\otimes\!\! I_n)(\tilde A\!\!+\!\!\tilde H\tilde C)(\mathbf{1}_N\!\!\otimes\!\! I_n)\!\!=\!\!A\!\!+\!\!\textstyle{\sum}_{i\!=\!1}^{N}\!\zeta_iH_iC_i,\\
&\!(T^{-1}\!\!\otimes\!I_n)(\mathcal{L}\!\otimes\! I_n)(T\!\otimes\!I_n)\!\!=\!\!(T^{-1}\!\!\mathcal{L}T)\!\otimes\! I_n\!\!=\!\!\left(\begin{smallmatrix}
    0 & \textbf{0}^{\top}_{N\!-\!1}\\
    \textbf{0}_{N\!-\!1} & \Delta
\end{smallmatrix}\right)\!\!\otimes\! I_n,
\end{aligned}
\vspace{-0.3cm}\end{equation*}
then 
\vspace{-0.05cm}\begin{equation}\label{eq:lindynamiceta}\begin{aligned}
  \dot \eta\!\!=\!\!
  \left(
\begin{smallmatrix}
\hat A & \hat B\\
\hat C & \Phi-\nu(\Delta\otimes I_n)\\
\end{smallmatrix}
\right)\!\!\eta.
\end{aligned}\vspace{-0.15cm}\end{equation}
where 
 $\hat A\!\!=\!\!A\!\!+\!\!\textstyle{\sum}_{i\!=\!1}^{N}\!\zeta_iH_iC_i$, $\Phi\!\!=\!\!(Y_2^{\top}\!\!\otimes\! I_n)(\tilde A\!\!+\!\!\tilde H\tilde C)(Y_1\!\!\otimes\! \!I_n)$,
$\hat B\!\!=\!\!(\zeta^{\top}\!\!\!\otimes\!\! I_n)\tilde H\tilde C(Y_1\!\!\otimes\!\! I_n)$,
$\hat C\!\!=\!\!(Y_2^{\top}\!\!\!\otimes\!\! I_n)\tilde H\tilde C(\mathbf{1}_N\!\!\otimes\!\! I_n)$.
%
%
Consider the Lyapunov candidate $V(\eta)\!\!=\!\!\eta^{\top}\!\!P\eta$ with $P\!\!=\!\!I_N\!\otimes\! P_a$ and $P_a$ defined by LMI \eqref{LMI_Global_Observ}, \eqref{LMI:DeltaP}. The derivative of $V$ along
\eqref{eq:lindynamiceta} is
\vspace{-0.3cm}\begin{equation*}\begin{aligned}
  \dot V\!\!=&\eta^{\top}\!\!\!
  \left(
\begin{smallmatrix}
P_a\hat A+\hat A^{\top}\!\!P_a & P_a\hat B+\hat C^{\top}\!\!(I_{N-1}\!\otimes P_a)\\
\hat B^{\top}\!\! P_a+(I_{N-1}\!\otimes P_a)\hat C_1 & (I_{N-1}\!\otimes P_a)\Phi+\Phi^{\top}\!\!(I_{N-1}\!\otimes P_a)\!-\!\nu(\Delta^{\top}\!\!+\!\Delta)\otimes P_a\\
\end{smallmatrix}
\right)\eta.
\end{aligned}\vspace{-0.3cm}\end{equation*}
Since $H_i\!\!=\!\!\frac{\bar{H}_i}{\zeta_i}$, then \eqref{LMI_Global_Observ} is fulfilled. Let \eqref{LMI:DeltaP} be fulfilled, then, using Schur Complement (\cite{Boyd_etall1994:Book}), we derive $\dot V\!\!\leq\!\!-2\rho V$ for a sufficiently large $\nu$. Thus we conclude error equation \eqref{eq:lin_error_dy} with $q\!\!=\!\!\textbf{0}$ is globally uniformly exponentially stable. Then, with respect to perturbation $q\!\!\in\!\!L^\infty\!(\R,\R^{n+p})$, the ISS property is derived straightforwardly. 
\vspace{-0.7cm}\end{pf}

Lemma \ref{thm_linear} implicitly proves the matrix inequality:
\vspace{-0.3cm}\begin{equation}\label{eq:lin_LMI}\begin{aligned}
  P\tilde A\!\!+\!\!\tilde A^\top\!\!\! P\!\!+\!\!P\tilde{H}\tilde{C}\!\!+\!\!\tilde{C}^\top\!\!\!\tilde{H}^\top\!\!\! P\!\!-\!\!\nu (\mathcal{L}\!\!+\!\!\mathcal{L}^\top\!\!)\!\otimes\! P_a\!\!+\!\!2\rho P\!\!\prec\!\!0,
\end{aligned}
\vspace{-0.3cm}\end{equation}
where $0\!\!\prec\!\!P\!\!=\!\!I_N\!\otimes\! P_a\!\!\in\!\!\R^{Nn\times Nn}$.
\vspace{-0.3cm}


Below we follow the idea of an upgrade of the linear observer to a homogeneous one to obtain our design scheme. For the centralized observer the upgrading procedure is developed in
\cite{wang2021generalized}. It suggests making gains of linear observers be scaled by a homogeneous term dependent on the norm of output vectors. In case of distributed observers
\eqref{eq_liner_observer} the gains $H_i$ and $\nu$ have to be modified in a similar way. 



\vspace{-0.4cm}
\section{Homogeneous Observer Design}\label{MR}
\vspace{-0.3cm}

The structure of the proposed distributed observer is similar to the linear distributed observer, which means there are no special constraints on the topology. 

\vspace{-0.3cm}
\subsection{Globally Homogeneous Distributed Observer}
\vspace{-0.3cm}
Let the observer be designed as follows
\vspace{-0.3cm}\begin{equation}\label{eq:fin_time_observer}
\dot{\hat{x}}_i \!\!=\!\!A\hat{x}_i\!+\!Bu\!+\!\gamma(t,\hat x_i)\!+\!g(|\omega_i|)H_i\omega_i\!+\!\nu \|\theta_i\|_\dn^\mu \theta_i,\;\;i\!\!=\!\!\overline{1,N}
\vspace{-0.3cm}\end{equation}
where $\omega_i$ and $\theta_i$ are defined as in \eqref{eq_liner_observer}, $g(|\omega_i|)\!\!=\!\!\exp(\mu(G_0\!\!+\!\!I_n)\ln\!|\omega_i|)$, $\|\cdot\|_{\dn}$ is a canonical homogeneous norm to be defined below, linear dilation $\dn$ is generated by $G_\dn\!\!=\!\!\mu G_0\!\!+\!\!I_n$ with $G_0$ defined in \eqref{eq:G_0}, we notice such a $G_\dn$ is anti-Hurwitz if $\mu\!\!>\!\!-1/{\tilde n}$, $\tilde n\!\!=\!\!\min\{k\}\!\!:\!\!(C,CA, \dots, CA^{k-1})$ has full rank (\cite{Polyakov2020:Book}). 
\vspace{-0.2cm}

For the proposed observer, the error equation \eqref{eq_aug} becomes
\vspace{-0.3cm}\begin{equation}\label{eq:error_eq_fin}\begin{aligned}
\dot{e} \!\!=&(\tilde A\!+\!\diag\{g(|\omega_i|)\}_{i=1}^{N}\tilde H \tilde C\!-\!\nu\diag\{\|\theta_i\|_\dn^\mu I_n\}_{i=1}^{N}(\mathcal{L}\!\otimes\! I_n))e\\&\!+\!\Gamma(t,\hat x, x)\!-\!\diag\{g(|\omega_i|)\}_{i=1}^{N}\tilde Hq_y\!-\!\tilde q_x,    
\end{aligned}
\vspace{-0.05cm}\end{equation}
with $\Gamma(t,\hat x, x)\!\!=\!\!(\Gamma_1^\top\!\!(t,\hat x_1, x),\dots,\Gamma_N^\top\!(t,\hat x_N, x))^\top\!\!\!$, $\hat x\!\!=\!\!(\hat x_1^\top,\dots,\hat x_N^\top)^\top\!\!\!$, $\Gamma_i(t,\hat x_i, x)\!\!=\!\!\gamma(t,\hat x_i)\!-\!\gamma(t,x)$, $\tilde q_x\!\!=\!\!\mathbf{1}_N\!\otimes\! q_x$.
\vspace{-0.3cm}
\begin{thm}\label{thm_fin_homo} {\it Let $q\!\!=\!\!\textbf{0}$ and $L_0\!\!=\!\!\textbf{0}$. Let $\mu\!\!>\!\!-1/{\tilde n}$. Let Assumption  \ref{assum2} be fulfilled and $H_i$, $\nu$ be selected as in Lemma \ref{thm_linear} with $Y$ and $P_a$ by solving LMI \eqref{LMI:DeltaP} together the following LMI
 \vspace{-0.3cm}\begin{equation}\label{LMI_Gd}\begin{aligned} 
&P_a\!\!\succ\!\!0,\;\;\;P_aG_\dn\!+\!G_\dn^\top P_a\!\succ\!0,\\
&  P_aA\!+\!A^\top\!\!P_a\!+\!YC\!+\!C^{\top}\!\!Y^{\top}\!\!\!+\!\rho (P_aG_\dn\!+\!G_\dn^\top P_a)\!\!\prec\!\!0,
 \end{aligned} 
 \vspace{-0.3cm}\end{equation} 
 with $\rho\!\!>\!\!0$.
Let the canonical homogeneous norm $\|\cdot\|_{\dn}$ be induced by the weighted Euclidean norm $\|\cdot\|_{P_a}$. Let   
 the nonlinear distributed observer \eqref{eq:fin_time_observer} operate under a strongly connected graph $\mathcal{G}$. Then the error equation \eqref{eq:error_eq_fin} is
globally uniformly finite-time (\textit{resp.,} nearly fixed-time)  stable provided $\mu\!\!<\!\!0$ (\textit{resp.,} $\mu\!\!>\!\!0$) is close enough to zero and $\gamma$ subjects to the following constraint: 
\vspace{-0.3cm}\begin{equation}\label{eq:con_gamma}
    \|\tilde\dn(\!-\!\ln\!\|\hat x\!\!-\!\!\tilde x \|_{\tilde\dn})(\hat \gamma(t,\hat x)\!\!-\!\!\tilde \gamma(t, x))\|_P\!\!\leq\!\!\tau \|\hat x\!\!-\!\!\tilde x\|^\mu_{\tilde\dn}, \;\forall t\!\!>\!\!0,\; \forall x_i,x\!\!\in\!\!\R^n,
\vspace{-0.05cm}\end{equation}
with $0\!\!<\!\!\tau\!\!<\!\!\tfrac{\rho}{3}$, $\tilde x\!\!=\!\!\textbf{1}_N\!\otimes\!x$, $\hat \gamma(t,\hat x)\!\!=\!\!( \gamma^\top\!\!(t,\hat x_1),\dots,\gamma^\top\!\!(t,\hat x_N))^\top$, $\tilde \gamma(t, x)\!\!=\!\!\textbf{1}_N\!\otimes\!\gamma(t,x)$,
 $\|\cdot\|_{\tilde\dn}$ is the canonical homogeneous norm induced by $\|\cdot\|_P$, $P\!\!=\!\!I_N\!\otimes\!P_a$, $\tilde \dn\!\!=\!\!I_N\!\otimes\!\dn$ is the linear dilation generated by $G_{\tilde\dn}\!\!=\!\!I_N\!\otimes\! G_\dn$.  The settling-time function is 
\vspace{-0.3cm}\begin{equation}\label{eq:settletime}
    T(e(0))\!\!\leq\!\!\tfrac{3}{\mu(3\tau\!-\!\rho)}\|e(0)\|_{\tilde \dn}^{-\mu}.
\vspace{-0.3cm}\end{equation}}
\end{thm}
\vspace{-0.8cm}
\begin{pf}
For $q\!\!=\!\!\textbf{0}$, the right-hand side of \eqref{eq:error_eq_fin} is continuous on $e\!\!\in\!\!\R^{Nn}$ for $\mu\!\!>\!\!-1/{\tilde n}$. Detailed analysis is in \ref{app:con}.
\vspace{-0.2cm}

Next, we prove $\|e\|_{\tilde\dn}$ is a Lyapunov function for the error equation. Indeed, consider $L_0\!\!=\!\!\textbf{0}$, equation \eqref{eq:G_0} gives  $A G_\dn\!\!=\!\!( \mu I_n\!\!+\!\!G_\dn)A$, $CG_\dn\!\!=\!\!C$, which imply $\dn(s)A\!\!=\!\!\exp(-\mu s)A\dn(s)$, $C\dn(s)\!\!=\!\!\exp(s)C$, $\forall s\!\!\in\!\!\R$. 
Using formula $\eqref{eq:procannorm}$, we derive
\vspace{-0.3cm}\begin{equation}\label{eq:fin_error}\begin{aligned}
  &\tfrac{d\|e\|_{\tilde\dn}}{dt}\!\!=\!\!\tfrac{\|e\|_{\tilde\dn}^{1+\mu}\!\!z^\top \!\!P(\tilde A\!+\!\tilde D(\tilde Cz)\tilde H\tilde C\!-\!\nu\tilde \Theta(z)(\mathcal{L}\!\otimes\! I_n))z\!+\!\|e\|_{\tilde\dn}z^\top \!\!P\tilde\dn(-\ln\!\|e\|_{\tilde\dn})\Gamma(t,\hat x, x)}{z^\top\!\! P(I_N\!\otimes\! G_\dn) z},\\ 
  &\tilde D(\tilde Cz)\!\!=\!\!\diag\{g(|C_iz_i|)\}_{i=1}^{N},
  \tilde \Theta(z)\!\!=\!\!\diag\{\|(\mathcal{L}_i\!\otimes\! I_n)z\|_\dn^\mu I_n\}_{i=1}^{N},\\
  &z\!\!=\!\!(z_1^{\top}\!\!,\dots,z_N^{\top})^{\top}\!\!\!=\!\!\tilde\dn(-\!\ln\!\|e\|_{\tilde\dn})e,\; P\!\!=\!\!I_N\!\otimes\! P_a,
\end{aligned}\vspace{-0.05cm}\end{equation}
the derivation on obtaining \eqref{eq:fin_error} is detailed in \ref{app:finite}. 
Equation \eqref{eq:fin_error} can be equivalently written as 
\vspace{-0.3cm}\begin{equation}\begin{aligned}
  \tfrac{d\|e\|_{\tilde\dn}}{dt}\!\!=\!\!\|e\|_{\tilde\dn}^{1+\mu}\!\!\!&\left(\tfrac{z^\top\!\! P(\tilde A+\tilde H\tilde C-\nu(\mathcal{L}\otimes I_n))z\!+\!\|e\|^{-\mu}_{\tilde\dn}\!\!z^\top \!\!P\tilde\dn(-\ln\!\|e\|_{\tilde\dn})\Gamma(t,\hat x, x)}{z^\top\!\! P(I_{N}\otimes G_\dn) z}\right.\\
  &+\left.\!\!\!\tfrac{z^\top\!\! P((\tilde D(\tilde Cz)-I_{Nn})\tilde H\tilde C)z+z^\top\!\! P(\nu(I_{Nn}-\tilde \Theta(z))(\mathcal{L}\otimes I_n))z}{z^\top\!\! P(I_{N}\otimes G_\dn) z}\right)\!\!\!.
\end{aligned}\vspace{-0.05cm}\end{equation}
Since $H_i$ is selected by solving LMI \eqref{LMI:DeltaP}, \eqref{LMI_Gd}, and $\nu$ is large enough, similar to inequality \eqref{eq:lin_LMI} implicitly obtained by Lemma \ref{thm_linear}, we have 
\vspace{-0.3cm}\begin{equation}\label{eq:sim_LMI}\begin{aligned}
  P\tilde A\!\!+\!\!\tilde A^\top\!\!\! P\!\!+\!\!P\tilde{H}\tilde{C}\!\!+\!\!\tilde{C}^\top\!\!\!\tilde{H}^\top\!\!\! P\!\!&-\!\!\nu \!(\mathcal{L}\!\!+\!\!\mathcal{L}^\top\!\!)\!\!\otimes\!\! P_a\\
  &+\!\!\rho I_N\!\!\otimes\!\! (P_aG_\dn\!\!+\!\!G_\dn^\top \!\! P_a)\!\!\prec\!\!0.
\end{aligned}
\vspace{-0.3cm}\end{equation}
Then we have 
\vspace{-0.3cm}\begin{equation}\begin{aligned}
  \tfrac{d\|e\|_{\tilde\dn}}{dt}\!\!\leq \!\!\|e\|_{\tilde\dn}^{1+\mu}&\!\!\left(-\rho\!+\!\tfrac{\|e\|^{-\mu}_{\tilde\dn}z^\top \!\!P\tilde\dn(-\ln\!\|e\|_{\tilde\dn})\Gamma(t,\hat x, x)}{z^\top\!\! P(I_{N}\otimes G_\dn) z}\right.\!\\&+\left.\!\!\!\tfrac{z^\top\!\! P((\tilde D(\tilde Cz)-I_{Nn})\tilde H\tilde C)z+z^\top\!\! P(\nu(I_{Nn}-\tilde \Theta(z))(\mathcal{L}\otimes I_n))z}{z^\top\!\! P(I_{N}\otimes G_\dn) z}\right)\!\!\!.
\end{aligned}\vspace{-0.05cm}\end{equation}
For $\mu\!\!=\!\!0$ LMI \eqref{LMI_Gd} reduce to LMI \eqref{LMI_Global_Observ}, so, by continuity, this LMI is feasible for $\mu$ close to zero.  Since $z\!\!=\!\!\tilde\dn(-\ln\!\|e\|_{\tilde\dn})e$ belongs to the unit sphere then  $\|(\tilde D(\tilde Cz)\!-\!I_{Nn})\tilde H\tilde C z\|_P\!\!\to\!\! 0$ and $\|(I_{Nn}\!-\!\tilde \Theta(z))\!(\mathcal{L}\!\otimes\! I_n)z\|_P\!\!\to\!\! 0$ as $\mu\!\!\to\!\! 0$ uniformly on $e\!\!\in\!\! \R^{Nn}$ (the detailed analysis is in \ref{app:muto0}). Hence, for $\mu$ being sufficiently close to zero, we have $\tfrac{z^\top\!\! P((\tilde D(\tilde Cz)\!-\!I_{Nn})\tilde H\tilde C)z}{z^\top\!\! P(I_{N}\otimes G_\dn) z}\!\!<\!\!\tfrac{\rho}{3}$, $\forall e\!\!\neq\!\! \textbf{0}$, and $\tfrac{z^\top\!\! P(\nu(I_{Nn}\!-\!\tilde \Theta(z))(\mathcal{L}\!\otimes\! I_n))z}{z^\top\!\! P(I_{N}\otimes G_\dn) z}\!\!<\!\!\tfrac{\rho}{3}$, $\forall e\!\!\neq\!\! \textbf{0}$. 
In addition, since \eqref{eq:con_gamma} holds, one derives $\tfrac{\|e\|^{-\mu}_{\tilde\dn}z^\top \!\!P\tilde\dn(-\ln\!\|e\|_{\tilde\dn})\Gamma(t,\hat x, x)}{z^\top\!\! P(I_{N}\otimes G_\dn) z}\!\!\leq\!\!\tau\!\!<\!\!\tfrac{\rho}{3}
$ for $\mu$ close to zero. 
Thus, we have $\tfrac{d\|e\|_{\tilde\dn}}{dt}\!\!< \!\!-(\tfrac{\rho}{3}\!\!-\!\!\tau)\|e\|^{1+\mu}_{\tilde\dn} $ and the finite-time (nearly fixed-time) stability of the error equation for $\mu\!\!<\!\!0$ $(\mu\!\!>\!\!0)$. The settling-time function is as \eqref{eq:settletime}.
\end{pf}

\vspace{-0.6cm}
\begin{remark}
\label{assum:gamma}
Let $\gamma(t,x)\!\!=\!\!E\psi(t,x)$, $E\!\!\in\!\!\R^{n\!\times\!m}$, $\psi\!\in\!C(\R\!\!\times\!\!\R^{n},\R^{m})$ and $G_\dn E\!\!=\!\!cE$, $c\!\!>\!\!0$. The latter implies $\dn(s)\gamma\!\!=\!\!\exp(cs)\gamma$, $\forall s\!\!\in\!\!\R$. 
In this case the condition 
\eqref{eq:con_gamma} becomes
\vspace{-0.3cm}\begin{equation}\label{eq:gamma_lip}
\|\hat \gamma(t,\hat x)\!\!-\!\!\tilde \gamma(t, x)\|_P\!\!\leq\!\!\tau \|\hat x\!\!-\!\!\tilde x\|^{\mu+c}_{\tilde\dn}, \;{\forall t\!\!>\!\!0,\; \forall x_i,x\!\!\in\!\!\R^n}.
\vspace{-0.3cm}\end{equation}
Taking into account inequality \eqref{eq:hom_norm_est} we conclude that 
the condition \eqref{eq:con_gamma} is fulfilled if $\gamma$ satisfies a H${\ddot{o}}$lder-like condition.
\end{remark}
\vspace{-0.3cm}
\subsection{Locally Homogeneous Distributed Observer}\label{sec:fix}
\vspace{-0.3cm}
Inspired by \cite{andrieu2008homogeneous,lopez2018finite}, we define a distributed observer 
combining homogeneous components of different degrees:
\vspace{-0.3cm}\begin{equation}\label{eq:fixed_time_observer}\begin{aligned}
  &\!\!\!\!\!\dot{\hat{x}}_i \!\!=\!\!A\hat{x}_i\!\!+\!\!Bu\!\!+\!\!\gamma(t,\hat x_i)\!\!+\!\!g(|\omega_i|)H_i\omega_i\!\!+\!\!\nu h(\theta_i) \theta_i,i\!\!=\!\!\overline{1,N}\\
  &\!\!\!\!\!g(|\omega_i|)\!\!=\!\! \tfrac{1}{2}\textstyle{\sum}_{k}\!\exp(\mu_k(G_0\!\!+\!\!I_n)\!\ln\!|\omega_i|),
  h(\theta_i)\!\!=\!\! \tfrac{1}{2}\textstyle{\sum}_{k}\!\|\theta_i\|_{\dn_k}^{\mu_k},
\end{aligned}
\vspace{-0.3cm}\end{equation}
where $k\!\!\in\!\!\{0\}\!\!\cup\!\!\{\infty\}$, $\omega_i$ and $\theta_i$ are defined in \eqref{eq_liner_observer}, 
dilation $\dn_k$ is generated by $G_{\dn_{k}}\!\!\!=\!\!\mu_k G_0\!+\!I_n$,  $G_0$ is defined in \eqref{eq:G_0}, $\mu_k\!\!>\!\!-1/{\tilde n}$ is a local homogeneity degree (\cite{andrieu2008homogeneous}) such that $G_{\dn_k}$ anti-Hurwitz, $\|\cdot\|_{\dn_k}$ is the canonical homogeneous norm to be defined below. With the proposed observer, the error equation \eqref{eq_aug} becomes
\vspace{-0.3cm}\begin{equation}\label{eq:error_eq_fix}\begin{aligned}
  \dot{e} \!\!=&\!\!\left(\tilde A\!\!+\!\!\diag\{g(|\omega_i|)\}_{i=1}^{N}\tilde H \tilde C\!\!-\!\! \nu\diag\{h(\theta_i)I_n\}_{i=1}^{N}(\mathcal{L}\!\!\otimes\! \!I_n)\right)\!\!e\\&\!\!+\!\!\Gamma(t,\hat x, x)\!\!-\!\!\diag\{g(|\omega_i|)\}_{i=1}^{N}\tilde Hq_y\!\!-\!\!\tilde q_x.
  \end{aligned}\vspace{-0.05cm}\end{equation}
If $(A,C)$ is observable, then the LMI
 \vspace{-0.3cm}\begin{equation}\label{LMI_Global_Observ_2}
 \begin{array}{l}
P_a\!\!\succ\!\!0,\;P_aA\!+\!A^\top\!\!P_a\!+\!YC\!+\!C^{\top}\!\!Y^{\top}\!\!\!+\!2\rho P_a\!\!\prec\!\!0,\\
P_aA\!+\!A^\top\!\!P_a\!+\!\frac{1}{2}\!YC\!+\!\frac{1}{2}C^{\top}\!\!Y^{\top}\!\!\!+\!\rho (P_aG_{\dn_k}\!+\!G^\top_{\dn_k} P_a)\!\!\prec\!\!0
  \end{array}
  \vspace{-0.3cm}\end{equation}
is feasible for some $\rho\!\!>\!\!0$. Indeed, taking $Y\!\!=\!\!-C^{\top}$ we conclude that
the feasibility of the LMI 
\vspace{-0.3cm}\[
P_a\!\!\succ\!\!0,\;\;\;P_aA\!+\!A^\top\!\!P_a\!-\!C^{\top}\!\!C\!\!\prec\!\!  0
\vspace{-0.3cm}\]
implies the feasibility of \eqref{LMI_Global_Observ_2} for a small $\rho\!\!>\!\!0$. The latter LMI is feasible due to the observability of $(A,C)$ (see, \textit{e.g.,} \cite{hespanha2018linear}).

\vspace{-0.3cm}
\begin{thm}\label{thm:bilimit_stab} {\it{Let $q\!\!=\!\!\textbf{0}$ and $L_0\!\!=\!\!\textbf{0}$. Let $\mu_k\!\!>\!\!-1/{\tilde n}$, $k\!\!\in\!\!\{0\}\!\!\cup\!\!\{\infty\}$. Let Assumption  \ref{assum2} be fulfilled and $H_i$, $\nu$ be selected as Lemma \ref{thm_linear} with $Y$ and $P_a$ being a solution of  the LMI \eqref{LMI_Global_Observ_2} and
\vspace{-0.3cm}\begin{equation}\label{eq:Gd_blimit}\begin{aligned}
\!\!\!\!P_a\!\!\succ\!\!0,\;(\Delta^\top\!\!\!+\!\Delta)\!\!\otimes\!\! P_a\!\!\succ\!\!0,\;P_a G_{\dn_k}\!\!+ \!G_{\dn_k}^\top\!\!P_a\!\!\succ\!\!0,\;k\!\!\in\!\!\{0\}\!\!\cup\!\!\{\infty\}.
\end{aligned}\vspace{-0.3cm}\end{equation}
Let the canonical homogeneous norm $\|\cdot\|_{\dn_k}$ be induced by the weighted Euclidean norm $\|\cdot\|_{P_a}$.} 
Let the nonlinear distributed observer \eqref{eq:fixed_time_observer} operate under a strongly connected graph $\mathcal{G}$. Then, the error equation \eqref{eq:error_eq_fix} is 
globally uniformly fixed-time stable provided the homogeneity degrees $\mu_0\!\!<\!\!0$, $\mu_\infty\!\!>\!\!0$ are close enough to zero and 
\vspace{-0.3cm}\begin{equation}\label{eq:con_gamma_fx}  
\!\!\!\!\!\!\!\left\{\begin{smallmatrix}
     \|\tilde\dn_0(-\ln\!\|\hat x-\tilde x \|_{\tilde\dn_0})(\hat \gamma(t,\hat x)-\tilde \gamma(t, x))\|_P\leq \tau\|\hat x-\tilde x\|^{\mu_0}_{\tilde\dn_0}\hfill, &\|\hat x-\tilde x\|_P< 1 \hfill \\
     \|\tilde\dn_{\infty}(-\ln\!\|\hat x-\tilde x \|_{\tilde\dn_\infty})(\hat \gamma(t,\hat x)-\tilde \gamma(t, x))\|_P\leq \tau \|\hat x-\tilde x\|^{\mu_{\infty}}_{\tilde\dn_{\infty}}\hfill, &\|\hat x-\tilde x\|_P> 1 \hfill \\     
     \|\hat \gamma(t,\hat x)-\tilde \gamma(t, x)\|_P\leq\tau\|\hat x-\tilde x\|_P,\hfill & \|\hat x-\tilde x\|_P\in(1/\ell,\ell)\hfill
\end{smallmatrix}\right.\!\!\!,
\vspace{-0.3cm}\end{equation}
 where $\hat{x}$ is defined in \eqref{eq:error_eq_fin}, $\tilde x$, $\hat{\gamma}$, $\tilde \gamma$ are defined in Theorem \ref{thm_fin_homo}, 
$\|\cdot\|_{\tilde\dn_k}$ is the canonical homogeneous norm induced by the weighted Euclidean norm $\|\cdot\|_{P}$, $P\!\!=\!\!I_N\!\otimes\!P_a$, $\tilde \dn_k$ is the linear dilation generated by $G_{\tilde\dn_k}\!\!\!\!=\!\!I_N\!\otimes\! G_{\dn_k}$, $k\!\!\in\!\!\{0\}\!\!\cup\!\!\{\infty\}$, $\ell\!\!>\!\!1$ is sufficiently large, and $0\!\!<\!\!\tau\!\!<\!\!\tfrac{\rho}{3}$,}
\end{thm}
\vspace{-0.6cm}
\begin{pf}
    Since \eqref{eq:G_0} holds, then $\dn_k(s)A\!\!=\!\!\exp(-\mu_k s)A\dn_k(s)$, $C\dn_k(s)\!\!=\!\!\exp(s)C$, $\forall s\!\!\in\!\!\R$. Let $\|e\|_{\tilde\dn_k}$ be the  Lyapunov candidate
    in $k$-limit ($k\!\!=\!\!0$ or $k\!\!=\!\!\infty$). Using formula $\eqref{eq:procannorm}$, we have the derivative of $\|e\|_{\tilde\dn_k}$ along $\dot{{e}} \!\!=\!\!f(e,\mathbf{0})$ being
\vspace{-0.3cm}\begin{equation}\label{eq:fix_error}\begin{aligned}
    \tfrac{d \|e\|_{\tdn_k}}{d t}\!\!=\!\!\tfrac{\|e\|_{\tdn_k}e^{\top}\!\tdn_k^{\top}\!(-\ln \|e\|_{\tdn_k}) P\tdn_k(-\ln \|e\|_{\tdn_k})\dot e}{e^{\top}\!\tdn_k^{\top}\!(-\ln \|e\|_{\tdn_k})PG_{\tdn_k}\tdn_k(-\ln \|e\|_{\tdn_k})e},
	\end{aligned}\vspace{-0.3cm}\end{equation}
 where
\vspace{-0.3cm}\begin{equation*}\begin{aligned}
\dot{e} \!\!=\!\!(\tilde A\!\!+\!\!\diag\{g(|&\omega_i|)\}_{i=1}^{N}\!\tilde H \tilde C)e\!\!
\\&-\!\!\tfrac{\nu}{2}\diag\{\textstyle{\sum}_k\!\|\theta_i\|_{\dn_k}^{\mu_k} I_n\}_{i=1}^{N}\!(\mathcal{L}\!\!\otimes\!\! I_n)e\!\!+\!\!\Gamma(t,\hat x, x).    
\end{aligned}
\vspace{-0.3cm}\end{equation*}
By calculation (details are in \ref{app:fix}), we derive
\vspace{-0.3cm}\begin{equation}\label{eq:dercnk}\begin{aligned}
  &\tfrac{d\|e\|_{\tilde\dn_k}}{dt}\!\!<\!\!\|e\|_{\tilde\dn_k}^{1\!+\!\mu_k}\!\times\\&\!\!\tfrac{\!\!\tfrac{\rho}{3}\!\!+\!z^\top\!\! P(\tilde A\!+\!\tfrac{1}{2}\!\tilde D_k\!\!(\tilde Cz)\tilde H\tilde C\!-\!  \tfrac{1}{2}\!\tilde \Theta_k\!\!(z)\!\!(\mathcal{L}\!\otimes\! I_n))z\!+\!\|e\|_{\tilde\dn_k}^{-\mu_k}\!\!z^\top\!\! P\tilde\dn_k\!\!(-\ln\!\|e\|_{\tilde\dn_k})\Gamma(t,\hat x, x)}{z^\top\!\! P(I_N\!\otimes\! G_{\dn_k}) z}, 
  \end{aligned}\vspace{-0.05cm}\end{equation}
  where $\tilde\dn_k\!\!=\!\!I_N\!\otimes\!\dn_k$, $P\!\!=\!\!I_N\!\otimes\! P_a$, $\tilde D_k(\tilde Cz)\!\!=\!\!\diag\{\exp(\mu_k(G_0\!\!+\!\!I_n)\!\ln\!|C_iz_i|)\}_{i=1}^N$, 
  $\tilde \Theta_k\!(z)\!\!=\!\!\diag\{\|(\mathcal{L}_i\!\!\otimes\! \!I_n)z\|_{\dn_k}^{\mu_k} I_n\}_{i=1}^{N}$, and
  $z\!\!=\!\!(z_1^{\top},\dots,z_N^{\top})^{\top}\!\!\!\!\!=\!\!\tilde\dn_k(-\ln\!\|e\|_{\tilde\dn_k})e$.
\eqref{eq:dercnk} is equivalent to
\vspace{-0.3cm}\begin{equation*}\begin{aligned}
  \tfrac{d\|e\|_{\tilde\dn_k}}{dt}\!\!<\!\!\|e\|_{\tilde\dn_k}^{1\!+\!\mu_k}\!\!\!&\left(\!\tfrac{z^\top\!\! P(\tilde A\!+\!\frac{1}{2}\tilde H\tilde C\!-\!\frac{\nu}{2}(\mathcal{L}\otimes I_n))z\!+\!\|e\|_{\tilde\dn_k}^{-\mu_k}\!\!z^\top\!\! P\tilde\dn_k(-\!\!\ln\!\|e\|_{\tilde\dn_k})\Gamma(t,\hat x, x)}{z^\top\!\! P(I_{N}\otimes G_{\dn_k}) z}\right.
  \\&\!\!+\!\!\left.\tfrac{\tfrac{\rho}{3}\!\!+\!\!\tfrac{1}{2}z^\top\!\! P(( \tilde D_{k}(\tilde Cz)\!-\!I_{Nn})\tilde H\tilde C\!+\!{\nu}(I_{Nn}\!-\!\tilde \Theta_{k}(z))(\mathcal{L}\otimes I_n))z}{z^\top\!\! P(I_{N}\otimes G_{\dn_k}) z}\right).
\end{aligned}\vspace{-0.3cm}\end{equation*}
Repeating the proof of Theorem $\ref{thm_fin_homo}$ we have
$\tfrac{z^\top\!\! P(\tilde D_k(\tilde Cz)\!-\!I_{Nn})\tilde H\tilde Cz}{z^\top\!\! P(I_{N}\otimes G_{\dn_k}) z}\!\!<\!\!\tfrac{\rho}{3}$, $\tfrac{\nu z^\top\!\! P(I_{Nn}\!-\tilde \Theta_k(z))(\mathcal{L}\!\otimes\! I_n)z}{z^\top\!\! P(I_{N}\otimes G_{\dn_k}) z}\!\!<\!\!\tfrac{\rho}{3}$, 
$\frac{\|e\|^{-\mu_k}_{\tilde\dn}z^\top \!\!P\tilde\dn_k(-\ln\!\|e\|_{\tilde\dn})\Gamma(t,\hat x, x)}{z^\top\!\! P(I_{N}\otimes G_{\dn_k}) z}\!\!\leq\!\tau\!\!<\!\!\tfrac{\rho}{3}
$, $\forall e\!\!\neq\!\! \textbf{0}$
provided  $\mu_k\!\!\to\!\!0$ and  $\gamma$ subjects to \eqref{eq:con_gamma_fx}. Thus, $\tfrac{d\|e\|_{\tilde\dn_k}}{dt}\!\!<\!\!(-\tfrac{\rho}{3}\!\!+\!\!\tau)\|e\|_{\tilde\dn_k}^{1\!+\!\mu_k}$ if 
$z^\top\!\! P(\tilde A\!\!+\!\!\tfrac{1}{2}\tilde H\tilde C\!\!-\!\!\tfrac{\nu}{2}(\mathcal{L}\!\!\otimes\!\! I_n))z\!\!<\!\!-\rho z^\top\!\! P(I_{N}\!\!\otimes\!\! G_{\dn_k}) z$, \textit{i.e.,} if matrix $\Pi\!\!=\!\!P(\tilde A\!\!+\!\!\tfrac{1}{2}\tilde H\tilde C\!\!-\!\!\tfrac{\nu}{2}(\mathcal{L}\!\!\otimes\! \!I_n)\!\!+\!\!\rho (I_{N}\!\!\otimes\!\! G_{\dn_k}))$ is  Hurwitz. 
Recall the proof of Lemma \ref{thm_linear}, the stability  of matrix $\Pi$ is guaranteed by $\nu$ and $H_i$ defined as Lemma \ref{thm_linear} with $Y$ and $P_a$ obtained by solving LMI \eqref{LMI_Global_Observ_2}, \eqref{eq:Gd_blimit}. 
In this case, error equation $\dot e\!\!=\!\!f(e,\textbf{0})$ is globally nearly fixed-time stable in $\infty$-limit (\cite{andrieu2008homogeneous}) for $\mu_\infty\!\!>\!\!0$, \textit{i.e.,} for $e(0)\!\in\!\R^{Nn}\!\!\Rightarrow\!\! e(t)\!\!\rightarrow\!\!{B}(R)$ for some $R\!\!>\!\!0$. Meanwhile, $\dot e\!\!=\!\!f(e,\textbf{0})$ is locally  finite-time stable in $0$-limit for $\mu_0\!\!<\!\!0$, \textit{i.e.,} $e(0)\!\!\in\!\!{B}(r)\!\! \Rightarrow\!\! e(t)\!\!\rightarrow\!\! \textbf{0}$ for some $r\!\!>\!\!0$.
\vspace{-0.3cm}


Let $V\!\!=\!\!\|e\|_P^2$ be the Lyapunov candidate whose derivative along $\dot e\!\!=\!\!f(e,\textbf{0})$ is given by
\vspace{-0.3cm}\begin{equation*}\begin{aligned}
  \tfrac{d\|e\|_P^2}{dt}\!\!=\!\!2e^\top\!\!\! P(\tilde A\!+\!\diag\{g(|&\omega_i|)\}_{i=1}^{N}\!\tilde H \tilde C)e\!+\!2e^\top\!\! P\Gamma(t,\hat x, x)\\&\!\!-\!\! 2\nu e^\top\!\!\! P\diag\{h(\theta_i)I_n\}_{i=1}^{N}\!(\mathcal{L}\!\!\otimes \!\!I_n)e.
\end{aligned}\vspace{-0.3cm}\end{equation*}
Let $\ell\!\!>\!\!1$ be an arbitrary number such that $\ell\!\!>\!\!\max\{1/r,R\}$. The closure of the set 
${B}(\ell)\backslash{B}(1/\ell)$ is a compact.
Repeating considerations from \ref{app:muto0}, we conclude $g(|\omega_i|)\!\!\to\!\! I_n$ and $h(\theta_i)\!\!\to\!\! 1$ as $\mu_k\!\!\to\!\! 0$ uniformly  on
$e\!\!\in\!\!{B}(\ell)\backslash{B}(1/\ell)$.  Then 
$\tfrac{d\|e\|_P^2}{dt}\!\!\to\!\! 2e^\top\!\! P(\tilde A\!\!+\!\!\tilde H\tilde C\!\!-\!\!\nu(\mathcal{L}\!\!\otimes\!\! I_n))e\!\!+\!\!2e^\top\!\! P\Gamma(t,\hat x, x)$
as  $\mu_k\!\!\to\!\! 0$ uniformly  on
$e\!\!\in\!\!{B}(\ell )\backslash{B}(1/\ell)$. In addition, \eqref{eq:con_gamma_fx} gives  $\|\Gamma(t,\hat x, x)\|_P\!\!\leq\!\!\tau\|e\|_P$ for $e\!\!\in\!\!{B}(\ell)\backslash{B}(1/\ell)$.  Therefore, $\frac{d\|e\|^2_P}{dt}\!\!<\!\!-2(1\!\!-\!\!\tau)\|e\|^2_P$ for all 
$e\!\!\in\!\!{B}(\ell)\backslash{B}(1/\ell)$ provided that $\mu_k$ is close enough to zero.
Since $\ell $ can be selected arbitrarily large then, taking into account the local finite-time stability around the origin and global nearly fixed-time stability in the infinity we complete the proof.
\end{pf}
\vspace{-0.6cm}
For $L_0\!\neq\!\textbf{0}$ the system matrix $A$ is not necessarily nilpotent, but the observer \eqref{eq:fixed_time_observer} still valid if $|L_0|$ is small enough.
\vspace{-0.3cm}\begin{cor}\label{coro:L0neq0}
 {\itshape Under conditions of Theorem \ref{thm:bilimit_stab} 
 the  distributed observer \eqref{eq:fixed_time_observer} is globally uniformly fixed-time stable for $L_0\!\!\neq\!\! \mathbf{0}$ provided $|L_0|$ is small enough.}
\end{cor}
\vspace{-0.6cm}\begin{pf}
Indeed, for $L_0\!\!\neq\!\! \mathbf{0}$ the error equation \eqref{eq_aug} becomes
\vspace{-0.3cm}\begin{equation}\label{eq:error_eq_fix_L0}\begin{aligned}
  \!\!\!\!\!\!\!\!\!\!\!\dot{e} \!\!=\!\!\tilde f(e,\textbf{0}&)\!\!=\!\!\left(\tilde A_0\!+\!\diag\{g(|\omega_i|)\}_{i=1}^{N}\!\tilde H \tilde C\right)e\!+\!\Gamma(t,\hat x, x)\\&-\!\left( \nu\diag\{h(\theta_i)I_n\}_{i=1}^{N}\!(\mathcal{L}\!\otimes \!I_n)\!+\!(I_N\!\otimes\!L_0C) \right)\!e.\!\!\!\!\!\!\!
  \end{aligned}\vspace{-0.3cm}\end{equation}
In this case, we have 
\vspace{-0.3cm}\begin{equation*}\begin{aligned}
|f(\tilde\dn_k(\ln\!\|e&\|_{\tilde\dn_k})e)\!-\! \tilde f(\tilde\dn_k(\ln\!\|e\|_{\tilde\dn_k})e)|\!\!=\\&
           \!\!|(I_N\!\otimes\!L_0C)\tilde\dn_k(\ln\!\|e\|_{\tilde\dn_k})e|\!\!\leq\!\!|(I_N\!\otimes\!L_0C)P^{-1/2}|,
\end{aligned}\vspace{-0.3cm}\end{equation*}
for any $e\!\!\in\!\!\R^{Nn} \backslash\{\mathbf{0}\}$,  $k\!\!\in\!\!\{0\}\!\!\cup\!\!\{\infty\}$. So, 
the proof of Theorem \ref{thm:bilimit_stab} and the estimates
$\tfrac{d\|e\|_{\tilde\dn_0}}{dt}\!\!<\!\!-(\tfrac{\rho}{3}\!\!-\!\!\tau)\|e\|_{\tilde\dn_0}^{1+\mu_0}$ for all $e\!\!\in\!\! B(r)$, $\tfrac{d\|e\|_{\tilde\dn_{\infty}}}{dt}\!\!<\!\!-(\tfrac{\rho}{3}\!\!-\!\!\tau)\|e\|_{\tilde\dn_\infty}^{1+\mu_\infty}$ for all $e\!\!\in\!\! \R^{Nn}\backslash B(R)$ and $\frac{d\|e\|^2_P}{dt}\!\!<\!\!-2(1\!\!-\!\!\tau)\|e\|^2_P$ for all $e\!\!\in\!\! B(\ell)\backslash B(1/\ell )$  with $\ell\!\!>\!\!\max\{1/r,R\}$ remain valid for a small enough $|L_0|$. 
\end{pf}

\vspace{-0.6cm}
\subsection{Robustness Analysis}
\vspace{-0.3cm}
Now we move to the perturbed case, \textit{i.e.,} $q\!\!\neq\!\!\textbf{0}$. The error equation for the finite- and fixed-time distributed observer are presented as \eqref{eq:error_eq_fin} and \eqref{eq:error_eq_fix}, respectively. 
\vspace{-0.3cm}\begin{Pro}\label{pro:finite}
{\it{   Let conditions of Theorem \ref{thm_fin_homo} hold. The error equation \eqref{eq:error_eq_fin} is ISS with respect to the bounded perturbation $q\!\!=\!\!(q^\top_x,q^\top_y)^\top\!\!\!\!\!\in\!\!L^\infty(\R, \R^{n+p})$.
}}
\end{Pro}
\vspace{-0.6cm}
\begin{pf}
The error equation \eqref{eq:error_eq_fin} can be rewritten as 
\vspace{-0.3cm}\begin{equation}\label{eq:error_eq_fin_robust}\begin{aligned}
\dot{e} \!\!=&\tilde Ae\!+\!\diag\{g(|\omega_i|)\}_{i=1}^{N}\!\tilde H (\tilde Ce\!\!-\!\!q_y)\!\!+\!\!\Gamma(t,\hat x, x)\!\!-\!\!\tilde q_x\\&\!\!-\!\!\nu\diag\{\|\theta_i\|_\dn^\mu I_n\}_{i=1}^{N}\!\!(\mathcal{L}\!\!\otimes\!\! I_n)e, \;\;\; \omega_i\!\!=\!\!C_ie_i\!\!-\!\!q_{y,i}.  
\end{aligned}
\vspace{-0.3cm}\end{equation}
We prove that $\|e\|_{\tilde\dn}$ is the ISS-Lyapunov function. Indeed, 
\vspace{-0.3cm}\begin{equation*}\begin{aligned}
  &\tfrac{d\|e\|_{\tilde\dn}}{dt}\!\!=\!\!\tfrac{\|e\|_{\tilde\dn}^{1+\mu}}{z^\top\!\!\! P(I_N\!\otimes\! G_\dn) z}\!\!\left(z^\top \!\!\!P\tilde Az\!\!+\!\!z^\top \!\!P\tilde D(\tilde Cz,\|e\|_{\tilde\dn}^{-1}\!q_y)\tilde H\epsilon\!\right.\\
  &\left.-\!\nu z^\top \!\!P\tilde \Theta(z)(\mathcal{L}\!\!\otimes\!\! I_n)z\!\!+\!\!\|e\|^{-\mu}_{\tilde\dn}\!z^\top \!\!P\tilde\dn(-\ln\!\|e\|_{\tilde\dn})(\Gamma(t,\hat x, x)\!\!-\!\!\tilde q_x)\!\!\!\right),\\ 
  &\tilde D(\tilde Cz,\|e\|_{\tilde\dn}^{-1}\!\!q_y)\!\!=\!\!\diag\{g(|\epsilon_i|)\}_{i=1}^{N},\;\;\epsilon_i\!\!=\!\!C_iz_i\!\!-\!\!\|e\|_{\tilde\dn}^{-1}\!\!q_{y,i},\\
  &\epsilon\!\!=\!\!(\epsilon_1^\top\!\!\!,\dots,\epsilon_p^\top)^\top\!\!\!=\!\!(\tilde Cz\!\!-\!\!\|e\|_{\tilde\dn}^{-1}\!q_y), \;\;P\!\!=\!\!I_N\!\otimes\! P_a,\\
  &\tilde \Theta(z)\!\!=\!\!\diag\{\|(\mathcal{L}_i\!\!\otimes\!\! I_n)z\|_\dn^\mu I_n\}_{i=1}^{N},\;\;
  z\!\!=\!\!\tilde\dn(-\!\ln\!\|e\|_{\tilde\dn})e,
\end{aligned}\vspace{-0.3cm}\end{equation*}
in which the derivation on obtaining the second term of the right-hand side of the latter equation is as \ref{app:rob_fin} while the others are identical to the proof of Theorem \ref{thm_fin_homo}.
Then,
\vspace{-0.3cm}\begin{equation}\label{eq:ly_de}\begin{aligned}
  &\tfrac{d\|e\|_{\tilde\dn}}{dt}\!\!=\!\!\tfrac{\|e\|_{\tilde\dn}^{1+\mu}}{z^\top\!\! P(I_N\!\otimes\! G_\dn) z}\!\left(z^\top \!\!P(\tilde A\!\!+\!\!\tilde H\tilde C\!\!-\!\!\nu (\mathcal{L}\!\!\otimes\!\! I_n))z\!\!-\!\!\|e\|_{\tilde\dn}^{-1}\!z^\top \!\!P\tilde Hq_y\right.\!\\
  &+\!\!z^\top \!\!P(\tilde D(\tilde Cz,\|e\|_{\tilde\dn}^{-1}\!q_y)\!\!-\!\!I_{Nn})\tilde H\epsilon\!\!+\!\!\nu z^\top \!\!P(I_{Nn}\!\!-\!\!\tilde \Theta(z))(\mathcal{L}\!\!\otimes\!\! I_n)z\!\!\\
  &\left.+\!\|e\|_{\tilde\dn}^{-\mu}\!\!z^\top \!\!\!P\tilde\dn(\!-\!\ln\!\|e\|_{\tilde\dn})(\Gamma(t,\hat x, x)\!\!-\!\!\tilde q_x)\!\!\right).
\end{aligned}\vspace{-0.05cm}\end{equation}
\vspace{-0.3cm}

In the right-hand side of \eqref{eq:ly_de}, firstly, we have 
$z^\top \!\!P(\tilde A\!\!+\!\!\tilde H\tilde C\!\!-\!\!\nu (\mathcal{L}\!\!\otimes\!\! I_n))z\!\!-\!\!\|e\|_{\tilde\dn}^{-1}\!z^\top \!\!P\tilde Hq_y\!\!<\!\!-\tfrac{8\rho}{9}$
if $z^\top \!\!P(\tilde A\!\!+\!\!\tilde H\tilde C\!\!-\!\!\nu (\mathcal{L}\!\!\otimes\!\! I_n))z\!\!<\!\!-\rho$ (guaranteed by Theorem \ref{thm_fin_homo} with $\mu$ close to zero)   and $\|e\|_{\tilde\dn}$ satisfies  $\|e\|_{\tilde\dn}\!\!>\!\!\tfrac{9\sqrt{p}| P^{{1}/{2}}\tilde H|}{\rho}\|q_y\|_{L_\infty}$.
Secondly, we have $z^\top\!\! P(\tilde D(\tilde Cz,\|e\|_{\tilde\dn}^{-1}\!q_y)\!\!-\!\!I_{Nn})\tilde H\epsilon\!\!<\!\!\frac{\rho}{9}$ provided $\|e\|_{\tilde\dn}\!\!>\!\!\tfrac{\sqrt{p}}{\pi\!\!-\!\!|\tilde CP^{-{1}/{2}}|}\|q_y\|_{L_\infty}$ (The detailed proof is as \ref{app:rob_fin_second}). 
Thirdly, let $\tilde\dn_{q}$ be a montone linear dilation generated by $G_{\tilde\dn_{q}}\!\!\!\!\!=\!\!I_N\!\!\otimes\!\!(\mu(G_0\!\!+\!\!I_n)\!\!+\!\!I_n)$ satisfying $P_q\!\!\succ\!\!0$, $P_qG_{\tilde\dn_{q}}\!\!\!+\!\!G_{\tilde\dn_{q}}^\top P_q\!\!\succ\!\!0$, let $\sigma_q\!\!:\!\!\R\!\!\to\!\!\R$ with $\sigma_q(\xi)\!\!=\!\!|P^{\frac{1}{2}}||P_{q}^{-\frac{1}{2}}|\xi^{\lambda_m}$ and let $\xi_M\!\!=\!\!\sigma_q^{-1}(\rho/9)$, $\sigma_q^{-1}$ is the inverse function of $\sigma_q$. We have $\|e\|_{\tilde\dn}^{-\mu}\!\!z^\top \!\!\!P\tilde\dn(-\!\ln\!\|e\|_{\tilde\dn})\tilde q_x\!\!<\!\!\frac{\rho}{9}$ provided $\|e\|_{\tilde\dn}\!\!>\!\!\Upsilon_M\|\tilde q_x\|_{\tilde\dn_{q}}$, $\Upsilon_M\!\!=\!\!\max\{1,\xi_M^{-1}\}$ (The detailed proof is as \ref{app:rob_fin_third}).
Finally,
together  $\nu z^\top \!\!P(I_{Nn}\!\!-\!\!\tilde \Theta(z))(\mathcal{L}\!\!\otimes\!\! I_n)z\!\!<\!\!\tfrac{\rho}{3}$ with $\mu\!\!\to\!\!0$, and $\|e\|^{-\mu}_{\tilde\dn}z^\top \!\!P\tilde\dn(-\ln\!\|e\|_{\tilde\dn})\Gamma(t,\hat x, x)\!\!\leq\!\!\tau\!\!<\!\!\tfrac{\rho}{3}
$ given by Theorem \ref{thm_fin_homo}, we conclude $\tfrac{d\|e\|_{\tilde\dn}}{dt}\!\!<\!\!-(\tfrac{\rho}{3}\!\!-\!\!\tau)\|e\|_{\tilde\dn}^{1+\mu}$. Based on the analysis above, $\|e\|_{\tilde\dn}$ is the ISS-Lyapunov function, and the error system \eqref{eq:error_eq_fin} is ISS with respect to $q\!\!=\!\!(q^\top_x,q^\top_y)^\top\!\!\!\in\!\!L^\infty(\R, \R^{n+p})$.
\vspace{-0.2cm}\end{pf}

\begin{Pro}\label{pro:fix}
{\it{    Let conditions of Theorem \ref{thm:bilimit_stab} and Corollary \ref{coro:L0neq0} hold. Then error equation \eqref{eq:error_eq_fix} is ISS with respect to the perturbation $q\!\in\!L^\infty(\R, \R^{n+p})$ for both $L_0\!\!=\!\!\textbf{0}$ and $L_0\!\!\neq\!\!\textbf{0}$.
}}
\end{Pro}\vspace{-0.6cm}
\begin{pf}
Let $L_0\!\!=\!\!\textbf{0}$.
the robustness analysis is finished if we prove the ISS-Lyapunov function for $e\!\!\in\!\! \R^{Nn}\backslash B(R)$ and $e\!\!\in\!\! B(\ell)\backslash B(1/\ell )$ exists,  $\ell\!\!>\!\!\max\{1/r,R\}$, $R\!\!>\!\!r\!\!>\!\!0$. To this end, in the $\infty$-limit, let $\|e\|_{\tdn_\infty}$ be the Lyapunov function, whose derivative along \eqref{eq:error_eq_fix} yields the following estimation 
\vspace{-0.3cm}\begin{equation}\label{eq:rob_fix_estimation}\begin{aligned}
  &\tfrac{d\|e\|_{\tilde\dn_\infty}}{dt}\!\!<\!\!\tfrac{\|e\|_{\tdn_\infty}^{1+\mu_\infty}}{z^\top\!\!\! P(I_N\!\otimes\! G_{\dn_\infty}) z}\!\times\\
  &\left(z^\top \!\!\!P\tilde Az\!\!+\!\!\tfrac{1}{2}z^\top \!\!\!P\tilde D_\infty(\tilde Cz,\|e\|_{\tilde\dn_\infty}^{-1}\!q_y)\tilde H(\tilde Cz\!\!-\!\!\|e\|_{\tilde\dn_\infty}^{-1}\!q_y)\!\!+\!\!\tfrac{\rho}{9}\!\!+\!\!\tfrac{\rho}{6}\!\right.\\
  &\!\!\left.-\!\tfrac{\nu}{2} z^\top \!\!P\tilde \Theta_\infty(z)(\mathcal{L}\!\!\otimes\!\! I_n)z\!\!+\!\!\|e\|^{-\mu_\infty}_{\tilde\dn_\infty}\!z^\top \!\!P\tilde\dn_\infty(-\!\!\ln\!\|e\|_{\tilde\dn_\infty})(\Gamma\!\!-\!\!\tilde q_x)\!\right),
  \end{aligned}\vspace{-0.05cm}\end{equation}
  where $\tilde\dn_\infty\!\!=\!\!I_N\!\otimes\!\dn_\infty$, $P\!\!=\!\!I_N\!\otimes\! P_a$, 
  $\tilde \Theta_\infty\!(z)\!\!=\!\!\diag\{\|(\mathcal{L}_i\!\!\otimes\! \!I_n)z\|_{\dn_\infty}^{\mu_\infty} I_n\}_{i=1}^{N}$, and
  $z\!\!=\!\!\tilde\dn_\infty(-\ln\!\|e\|_{\tilde\dn_\infty})e$. Detailed derivation on obtaining \eqref{eq:rob_fix_estimation} is in \ref{app:rob_fix}.
  Then
\vspace{-0.3cm}\begin{equation*}\begin{aligned}
  &\tfrac{d\|e\|_{\tilde\dn_\infty}}{dt}\!\!<\\
  &\!\!\tfrac{\|e\|_{\tdn_\infty}^{1+\mu_\infty}}{z^\top\!\! P(I_N\!\otimes\! G_{\dn_\infty}) z}\!\!\left(z^\top \!\!P(\tilde A\!\!+\!\!\tfrac{1}{2}\tilde H\tilde C\!\!-\!\!\tfrac{\nu}{2}(\mathcal{L}\!\otimes\!I_n))z\!\!-\!\!\tfrac{1}{2}\|e\|_{\tilde\dn_\infty}^{-1}z^\top \!\!P\tilde Hq_y\right.\\
  &\left.\!\!+\!\!\tfrac{1}{2}z^\top \!\!P(\tilde D_\infty(\tilde Cz,\|e\|_{\tilde\dn_\infty}^{-1}\!q_y)\!\!-\!\!I_{Nn})\tilde H(\tilde Cz\!\!-\!\!\|e\|_{\tilde\dn_\infty}^{-1}\!q_y)\!\!+\!\!\tfrac{\rho}{9}\!\!+\!\!\tfrac{\rho}{6}\!\right.\\
  &\left.\!\!\!+\!\tfrac{\nu}{2} z^\top \!\!\!P(I_{Nn}\!\!-\!\!\tilde \Theta_\infty\!(z)\!)\!(\mathcal{L}\!\!\otimes\!\! I_n)z\!\!+\!\!\|e\|^{-\mu_\infty}_{\tilde\dn_\infty}\!z^\top \!\!P\tilde\dn_\infty\!\!(-\!\!\ln\!\|e\|_{\tilde\dn_\infty}\!\!)\!(\Gamma\!\!-\!\!\tilde q_x)\!\!\!\right).
  \end{aligned}\vspace{-0.3cm}\end{equation*}
Repeating the consideration of the proof of Proposition \ref{pro:finite}, 
we have $z^\top \!\!P(\tilde A\!\!+\!\!\tfrac{1}{2}\tilde H\tilde C\!\!-\!\!\tfrac{\nu}{2}(\mathcal{L}\!\!\otimes\!\! I_n))z\!\!-\!\!\tfrac{1}{2}\|e\|_{\tilde\dn_\infty}^{-1}\!z^\top \!\!P\tilde Hq_y\!\!<\!\!-\tfrac{17\rho}{18}$
if $z^\top \!\!P(\tilde A\!\!+\!\!\tfrac{1}{2}\tilde H\tilde C\!\!-\!\!\tfrac{\nu}{2} (\mathcal{L}\!\!\otimes\!\! I_n))z\!\!<\!\!-\rho$ (as Theorem \ref{thm:bilimit_stab} with $\mu_\infty$ close to zero) and $\|e\|_{\tilde\dn_\infty}$ satisfies  $\|e\|_{\tilde\dn_\infty}\!\!>\!\!\tfrac{9\sqrt{p}| P^{{1}/{2}}\tilde H|}{\rho}\|q_y\|_{L_\infty}$. In addition, similar to the proof of Proposition \ref{pro:finite}, $\tfrac{1}{2}z^\top\!\!\! P(\tilde D_\infty(\tilde Cz,\|e\|_{\tdn_\infty}^{-1}\!q_y)\!\!-\!\!I_{Nn})\tilde H(\tilde Cz\!\!-\!\!\|e\|_{\tdn_\infty}^{-1}\!q_y)\!\!<\!\!\frac{\rho}{18}$ provided $\|e\|_{\tilde\dn_\infty}\!\!>\!\!\tfrac{\sqrt{p}}{\pi\!\!-\!\!|\tilde CP^{-{1}/{2}}|}\|q_y\|_{L_\infty}$, $\pi\!\!>\!\!{|\tilde CP^{-{1}/{2}}|}$. Besides,
$\|e\|_{\tilde\dn_\infty}^{-\mu_\infty}\!\!z^\top \!\!\!P\tilde\dn_\infty(-\!\ln\!\|e\|_{\tilde\dn_\infty})\tilde q_x\!\!<\!\!\frac{\rho}{9}$ provided $\|e\|_{\tilde\dn_\infty}\!\!>\!\!\Upsilon_M\|\tilde q_x\|_{\tilde\dn_{q}}$, $\tilde\dn_{q}$ is a linear dilation generated by $G_{\tilde\dn_{q}}\!\!=\!\!I_N\!\otimes\!(\mu_\infty(G_0\!\!+\!\!I_n)\!\!+\!\!I_n)$, and $\Upsilon_M$ is defined in the proof of Proposition \ref{pro:finite}. Together  $\frac{\nu}{2} z^\top \!\!P(I_{Nn}\!\!-\!\!\tilde \Theta_\infty(z))(\mathcal{L}\!\!\otimes\!\! I_n)z\!\!<\!\!\tfrac{\rho}{6}$ with $\mu_\infty\!\!\to\!\!0$, and $\|e\|^{-\mu_\infty}_{\tilde\dn_\infty}z^\top \!\!P\tilde\dn_\infty(-\ln\!\|e\|_{\tilde\dn_\infty})\Gamma(t,\hat x, x)\!\!\leq\!\!\tau\!\!<\!\!\tfrac{\rho}{3}
$, we conclude
$\tfrac{d\|e\|_{\tdn_\infty}}{dt}\!\!<\!\!-(\tfrac{\rho}{3}\!\!-\!\!\tau)\|e\|_{\tilde\dn_\infty}^{1+\mu_\infty}$.
\vspace{-0.3cm}

On the other hand, for $e\!\!\in\!\!{B}(\ell)\backslash{B}(1/\ell)$, $\ell\!\!>\!\!\max\{1/r,R\}$, since $\mu_k$ is close enough to zero, the canonical homogeneous norm $\|\cdot\|_{\tilde\dn}$ reduces to $\|\cdot\|_{P}$.
In this case, we prove $\|e\|_P^2$ is the ISS-Lyapunov function, similar to the proof of Theorem \ref{thm:bilimit_stab}, we have
$\tfrac{d\|e\|_P^2}{dt}\!\!\to\!\! 2e^\top\!\! P(\tilde A\!+\!\tilde H\tilde C\!-\!\nu(\mathcal{L}\!\otimes\! I_n))e\!+\!2e^\top\!\! P\Gamma(t,\hat x, x)\!\!-\!\!2e^\top\!\! P\tilde Hq_y\!\!-\!\!2e^\top\!\! P\tilde q_x$
as  $\mu_k\!\!\to\!\! 0$ uniformly  on
$e\!\!\in\!\!{B}(\ell)\backslash{B}(1/\ell)$.
Since Theorem \ref{thm:bilimit_stab} holds, we have $e^\top\!\! P(\tilde A\!\!+\!\!\tilde H\tilde C\!\!-\!\!\nu(\mathcal{L}\!\!\otimes\!\! I_n))e\!\!<\!\!-\rho\|e\|^2_P$ and $e^\top\!\! P\Gamma(t,\hat x, x)\!\!\leq\!\!\!\tau\|e\|_P^2$. Moreover,
$e^\top\!\! P\tilde Hq_y\!\!<\!\!\tfrac{\rho}{3}\|e\|^2_P$ provided 
$\|e\|_P\!\!>\!\!\tfrac{3\sqrt{p}| P^{{1}/{2}}\tilde H|}{\rho}\|q_y\|_{L_\infty}$ and $e^\top\!\! P\tilde q_x\!\!<\!\!\tfrac{\rho}{3}\|e\|^2_P$ if $\|e\|_P\!\!>\!\!\tfrac{3\sqrt{Nn}|P^{{1}/{2}}|}{\rho}\| q_x\|_{L_\infty}$.
Therefore, we have
$\tfrac{d\|e\|_P^2}{dt}\!\!<\!\!-2(\tfrac{\rho}{3}\!\!-\!\!\tau)\|e\|^2_P$ provided $\|e\|_P\!\!>\!\!\max\{\tfrac{3\sqrt{p}| P^{{1}/{2}}\tilde H|}{\rho}\|q_y\|_{L_\infty},\tfrac{3\sqrt{Nn}|P^{{1}/{2}}|}{\rho}\| q_x\|_{L_\infty}\}$. Thus we can conclude the ISS stability of the error equation \eqref{eq:error_eq_fix} for $L_0\!\!=\!\!\textbf{0}$. On the other hand, for $L_0\!\!\neq\!\!\textbf{0}$, by repeating the consideration in Corollary \ref{coro:L0neq0}, the ISS of \eqref{eq:error_eq_fix} can be obtained straightforwardly.
   \vspace{-0.6cm}
\end{pf}

\vspace{-0.3cm}
\section{Simulation Results}\label{SR}
\vspace{-0.3cm}

Let the system matrices of the  plant \eqref{eq:LTI_plant} be
\vspace{-0.3cm}
\begin{equation*}\begin{aligned}
    A\!\!=\!\!\left(
      \begin{smallmatrix}
        0 & 1 & 0 \\
        0 & 0 & 1 \\
        0 & 0 & 0 \\
      \end{smallmatrix}
    \right),\;\; B\!\!=\!\!\left(
                 \begin{smallmatrix}
                   0 \\
                   0 \\
                   1 \\
                 \end{smallmatrix}
               \right).
  \end{aligned}
\vspace{-0.3cm}
\end{equation*}
Let $\gamma(x)\!\!=\!\!0.02(0\;1\;0)^{\top}\!\!|x|^{0.1}$.
Let the output matrices corresponding to $y_i,\;i\!\!=\!\!\overline{1,3}$ be $C_1\!\!=\!\!\left(
      \begin{smallmatrix}
      0 & 0 & 2 \\
      0 & 0 & 2
      \end{smallmatrix}
    \right)$,  $C_2\!\!=\!\!\left(
      \begin{smallmatrix}
      0 & 0 & 3
      \end{smallmatrix}
    \right)$,  $C_3\!\!=\!\!\left(
      \begin{smallmatrix}
      0 & 1 & 0\\
      3 & 2 & 2
      \end{smallmatrix}
     \right)$,
then we have $C\!\!=\!\!(C_1^\top,C_2^\top,C_3^\top)^\top$, and identity \eqref{eq:G_0} gives $G_0\!\!=\!\!\textbf{0}$, the latter implies $G_\dn\!\!=\!\!\mu G_0\!\!+\!\!I_n\!\!=\!\!I_n$. 
\vspace{-0.2cm}

\begin{figure}[htbp]
  \centering
  \includegraphics[width=2cm]{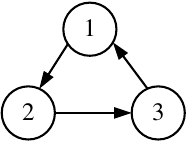}\\
  \vspace{-0.3cm}
  \caption{The communication graph of the distributed observer.}\label{fig:graph}
\vspace{-0.2cm}
\end{figure}

The distributed observer is composed of $3$ observers, whose communication graph is as Fig. \ref{fig:graph} shows.
The left zero-eigenvector of the associated Laplacian matrix is $\zeta\!\!=\!\!\frac{1}{3}\textbf{1}_3$. The initial states of observers are assigned to be zeros and the initial state of the plant is $x(0)\!\!=\!\! (-1.0\;0.0\;1.0)^{\top}\!$. 
\vspace{-0.15cm}

\vspace{-0.4cm}
\subsection{On the Robust Finite-time Distributed Observer}
\begin{figure}[htbp]
  \centering
  \includegraphics[width=7.5cm]{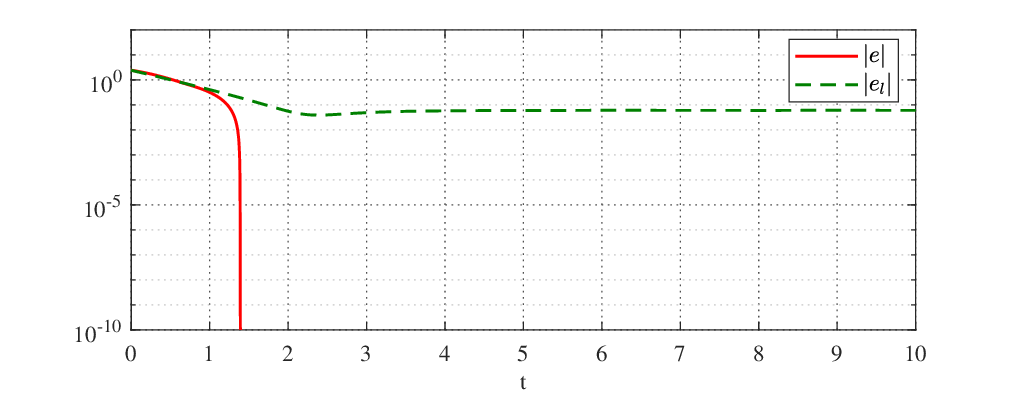}\\
  \vspace{-0.3cm}
  \caption{The trajectory of $|e|$ and $|e_l|$, with $q\!\!=\!\!\mathbf{0}$, by employing the finite-time distributed observer and the linear distributed observer, respectively.}\label{fig:normefin}
\vspace{-0.3cm}
\end{figure}

\vspace{-0.3cm}
Let $\mu\!\!=\!\!-0.65$ and $\nu\!\!=\!\!10$.
Let $\rho\!\!=\!\!1$, solving LMI \eqref{LMI_Global_Observ}, \eqref{LMI:DeltaP}, \eqref{LMI_Gd}  has 
\vspace{-0.3cm}$$\bar H_1\!\!=\!\!\left(\begin{smallmatrix}
    3.15 & -0.00\\
    -1.50 & -0.00\\
    -4.71 & -0.00\\
\end{smallmatrix}\right),\;\;
\bar H_2\!\!=\!\!\left(\begin{smallmatrix}
    -0.00 \\
    -0.00 \\
    0.00 \\
\end{smallmatrix}\right),\;\;
\bar H_3\!\!=\!\!\left(\begin{smallmatrix}
    3.30 & -3.15\\
    -9.37 & -0.00\\
    -0.00 & -0.00\\
\end{smallmatrix}\right).\vspace{-0.3cm}$$
Then, we can obtain $H_i\!\!=\!\!3\bar H_i$, $i\!\!=\!\!\overline{1,3}$ since $\zeta_i\!\!=\!\!\frac{1}{3}$, $i\!\!=\!\!\overline{1,3}$. 
%
Let the iteration step $h\!\!=\!\!0.001\mathrm{s}$ with iteration number $N\!\!=\!\!10000$, and  simulation is performed on the MATLAB platform.
The estimation error is $e\!\!=\!\!(e_1^{\top},e_2^{\top}, e_3^{\top})^{\top}\!$, where $e_i\!\!=\!\!\hat x_i\!-\!x$. In parallel, as a comparison, the estimation error for the classical linear distributed observer is defined as $e_l\!\!=\!\!(e_{l,1}^{\top},e_{l,2}^{\top}, e_{l,3}^{\top})^{\top}\!$. The comparison trajectory of $|e|$ and $|e_l|$ are as Fig. \ref{fig:normefin} with the observed plant initialized at $x(0)$. 
\vspace{-0.15cm}

\begin{figure}[htbp]
  \centering
  \includegraphics[width=7.5cm]{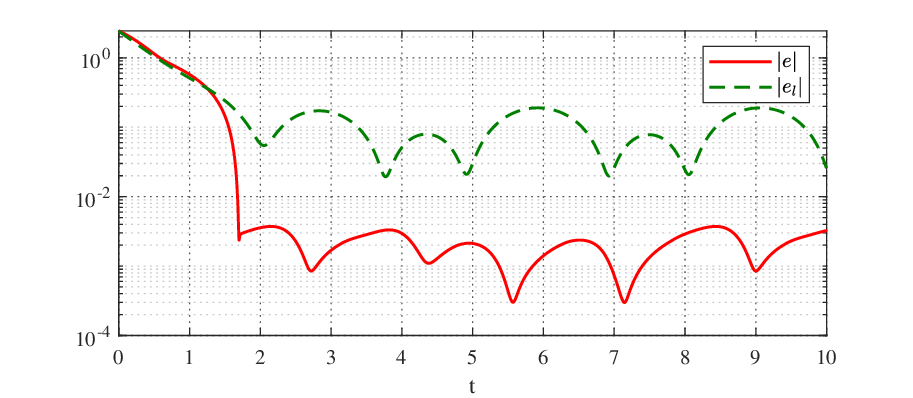}\\
  \vspace{-0.3cm}
  \caption{The trajectory of $|e|$ and $|e_l|$, with $q\!\!\neq\!\!\mathbf{0}$, by employing the finite-time distributed observer and the linear distributed observer, respectively.}\label{fig:normefinper}
\vspace{-0.2cm}
\end{figure}

 Then we concentrate on the robustness of the finite-time distributed observer. Let
 $q_x\!\!=\!\!0.1 (0\; 0\; \sin(2t))^\top$,
     $q_{y,1}\!\!=\!\!0.001(\sin(2t)\; \cos(0.5t))^\top$, $q_{y,2}\!\!=\!\!0.001\cos(t)$, 
    $q_{y,3}\!\!=\!\!0.001(\cos(2t)\; \sin(t))^\top$. The comparison trajectory of $|e|$ and $|e_l|$ is as Fig. \ref{fig:normefinper}. It is obvious that with bounded uncertainties in plant states as well as output measurements, the proposed finite-time observer can make estimation errors converge to a ball centering the origin with a decreased radius compared to the linear observer. 
\vspace{-0.15cm}

\vspace{-0.4cm}
\subsection{On the Robust Fixed-time Distributed Observer}
\begin{figure}[htbp]
  \centering
  \includegraphics[width=9cm]{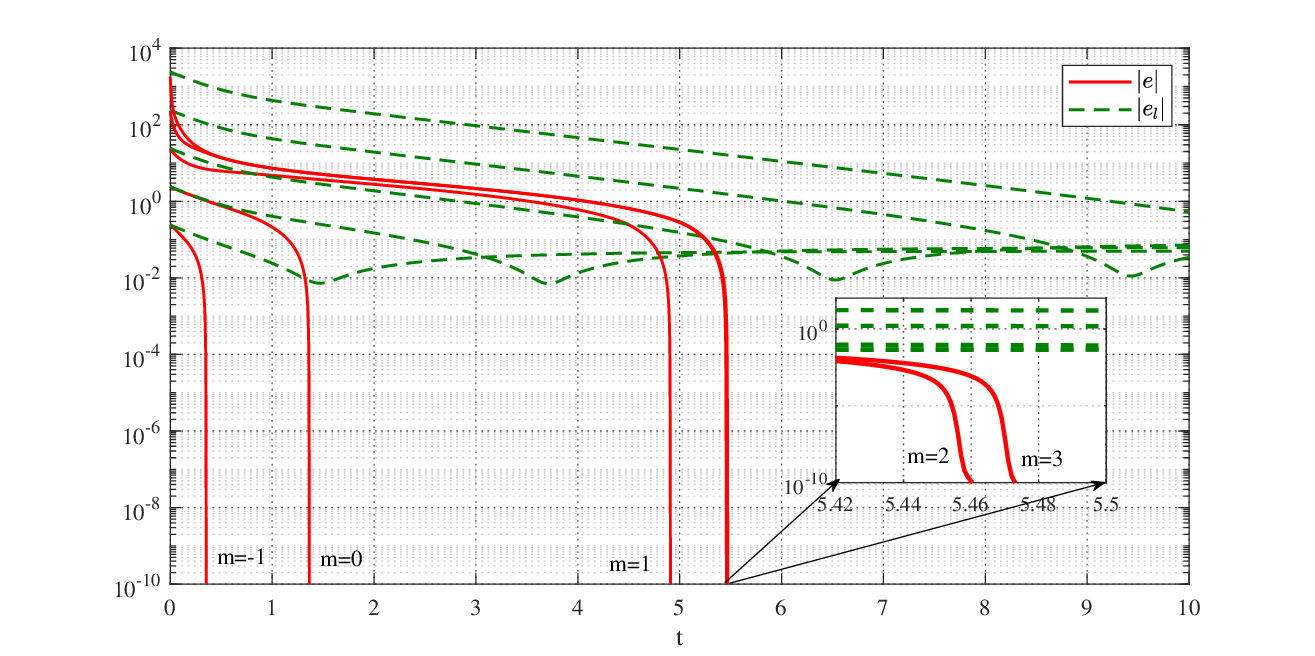}\\
  \vspace{-0.3cm}
  \caption{The trajectory of $|e|$ and $|e_l|$ by employing the fixed-time distributed observer and the linear distributed observer, respectively; with $q\!\!=\!\!\mathbf{0}$ and $m\!\!\in\!\!\{-1,0,1,2,3\}$.}\label{fig:normefix}
\vspace{-0.1cm}
 \end{figure}

\vspace{-0.3cm}
Let $\mu_0\!\!=\!\!-0.65$, $\mu_\infty\!\!=\!\!0.65$ and $\nu\!\!=\!\!10$.
Then, by solving LMI  \eqref{LMI_Global_Observ_2}, \eqref{eq:Gd_blimit} with $\rho\!\!=\!\!1$, one has 
\vspace{-0.3cm}$$\bar H_1\!\!=\!\!\left(\begin{smallmatrix}
     3.63  &  -0.00\\
    -2.60  &  -0.00\\
    -5.44  &  -0.00\\
\end{smallmatrix}\right),\;\;
\bar H_2\!\!=\!\!\left(\begin{smallmatrix}
    0.00 \\
    0.00 \\
    0.00 \\
\end{smallmatrix}\right),\;\;
\bar H_3\!\!=\!\!\left(\begin{smallmatrix}
    1.69  &  -3.69\\
    -10.33  & -0.23\\
    -0.68  & -0.02\\
\end{smallmatrix}\right).\vspace{-0.3cm}$$
Then $H_i\!\!=\!\!3\bar H_i$, $i\!\!=\!\!\overline{1,3}$. 
%
Let the iteration step $h\!\!=\!\!0.001\mathrm{s}$ with iteration number $N\!\!=\!\!10000$. The comparison trajectory of $|e|$ and $|e_l|$ is as Fig. \ref{fig:normefix} with the observed plant initializes at $10^mx(0)$, $m\!\!\in\!\!\{-1,0,1,2,3\}$. Fig. \ref{fig:normefix} confirms that the fixed-time design shows lower sensitivity in convergence time concerning the initial states of the observed plant.
\vspace{-0.3cm}

\begin{figure}[htbp]
  \centering
  \includegraphics[width=7cm]{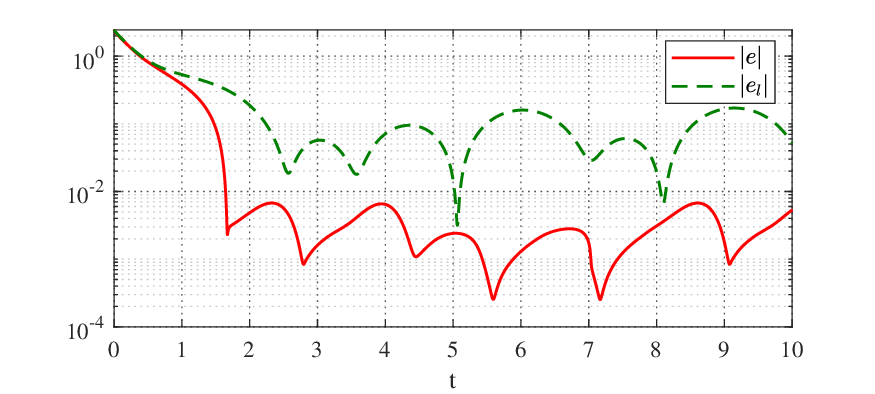}\\
  \vspace{-0.3cm}
  \caption{The trajectory of $|e|$ and $|e_l|$, with $q\!\!\neq\!\!\mathbf{0}$, by employing the fixed-time distributed observer and the linear distributed observer, respectively.}\label{fig:normefixper}
\vspace{-0.1cm}
\end{figure}

By employing the same perturbation defined in the finite-time simulation example. For the fixed-time distributed observer, the comparison trajectory of $|e|$ and $|e_l|$ is as Fig. \ref{fig:normefixper} with the observed plant initializes at $x(0)$. Compared to the linear case, 
 the error of the proposed fixed-time distributed observer converges to a smaller neighborhood of the origin. 
\vspace{-0.15cm}

\vspace{-0.0cm}
\section{Conclusion}
\vspace{-0.4cm}
In this paper, we tackle the problem of robust finite/fixed-time cooperative state estimation for a class of nonlinear systems by proposing a scheme based on distributed nonlinear observers, which is obtained by adopting the idea of upgrading the classical linear distributed observers to a generalized homogeneous one. By proper parameter tuning with LMIs, the proposed scheme offers the capability of achieving finite/fixed-time state reconstruction as well as admitting the nonlinearity yields H$\Ddot{\mathrm{o}}$lder conditions. The robustness with respect to bounded perturbations from both plant states and output measurements is also investigated.
\vspace{-0.4cm}
\begin{ack}
\vspace{-0.4cm}
This work was supported in part by the China Scholarship Council
(CSC) under Grant 202106160022.  
\end{ack}

\vspace{-0.4cm}

\appendix
\renewcommand{\thesection}{Appendix \Alph{section}}
\vspace{-0.0cm}
\section{}\label{app:con}
\vspace{-0.3cm}
In \eqref{eq:error_eq_fin}, for  $q\!\!=\!\!\textbf{0}$, the only possible discontinuous point is $e\!\!=\!\!\textbf{0}$. 
Notice that
\vspace{-0.3cm}
$$\diag\{g(|\omega_i|)\}_{i=1}^{N}\tilde H \tilde Ce\!\!=\!\!(\dots, (g(|\omega_i|)H_i\omega_i)^\top,\dots)^\top,
\vspace{-0.3cm}$$
where $\omega_i\!\!=\!\!C_ie_i$ and
\vspace{-0.3cm}\begin{equation*}\begin{aligned}
g(|\omega_i|)H_i\omega_i\!\!=\!\!\exp((I_n\!\!+\!\!\mu(G_0\!\!+\!\!I_n))\ln|\omega_i|)H_i\tfrac{\omega_i}{|\omega_i|},\;\;i\!\!=\!\!\overline{1,N}.
\end{aligned}
\vspace{-0.3cm}\end{equation*}
Since $\mu\!\!>\!\!-1/{\tilde n}$ and $G_0$ is selected by \eqref{eq:G_0}, then $I_n\!+\!\mu(G_0\!+\!I_n)$ is anti-Hurwitz \cite{Polyakov2020:Book}, \cite{zimenko2020robust} and $g(|\omega_i|)H_i\omega_i\!\!\to\!\!\textbf{0}$ as $|\omega_i|\!\!\to\!\!{0}$, $\forall i\!\!=\!\!\overline{1,N}$, the latter implies 
$\diag\{g(|\omega_i|)\}\tilde H\omega\!\!\to\!\!\textbf{0}$ as $|\omega|\!\!\to\!\!\textbf{0}$, so we have the function $e \!\!\to\!\! \diag\{g(|\omega_i|)\}_{i=1}^{N}\tilde H \tilde Ce$ is continuous at $e\!\!=\!\!\textbf{0}$.
\vspace{-0.2cm}

In addition, notice that 
\vspace{-0.3cm}\begin{equation*}\begin{aligned}
\nu\diag\{\|\theta_i\|_\dn^\mu I_n\}_{i=1}^{N}\!(\mathcal{L}\!\otimes\! I_n)e\!\!=\!\!(\dots, \nu\|\theta_i\|_\dn^\mu\theta_i^\top,\dots)^\top,
\end{aligned}\vspace{-0.3cm}\end{equation*}
with $\theta_i\!\!=\!\!(\mathcal{L}_i\!\otimes\!I_n)e$, $v\!\!>\!\!0$ and \vspace{-0.3cm}\begin{equation*}\begin{aligned}
\|\theta_i\|_\dn^\mu\theta_i\!\!=\!\!\exp((\mu (I_n\!+\!G_0)\!+\!I_n)\ln\!\|\theta_i\|_\dn)\dn(-\ln\!\|\theta_i\|_\dn)\theta_i.
\end{aligned}\vspace{-0.3cm}\end{equation*}
Using Cauchy-Schwarz inequality,
\vspace{-0.3cm}\begin{equation*}\begin{aligned}
0\!\!\leq\!\!\|\theta_i\|_\dn^\mu|\theta_i|\!\!\leq\!\!|\exp((\mu (I_n\!+\!G_0)\!+\!I_n)\ln\!\|\theta_i\|_\dn)||P_a^{-{1}/{2}}|.
\end{aligned}\vspace{-0.3cm}\end{equation*}
In the latter inequality, matrix $\mu(I_n\!+\!G_0)\!+\!I_n$ is anti-Hurwitz since $\mu\!\!>\!\!-1/{\tilde n}$ and $G_0$ is selected by \eqref{eq:G_0}.
Then, 
\vspace{-0.3cm}\begin{equation*}
    \begin{aligned}
|\exp((\mu (I_n\!+\!G_0)\!+\!I_n)\ln\!\|\theta_i\|_\dn)||P_a^{-{1}/{2}}|\!\!\to\!\!0,\;as\;{\|\theta_i\|_\dn\!\to\!{0}}.
    \end{aligned}
\vspace{-0.3cm}\end{equation*}
Therefore, we say
$\|\theta_i\|_\dn^\mu\theta_i\!\!\to\!\!\textbf{0}$ as $\|\theta_i\|_\dn\!\!\to\!\!{0}$, $\forall i\!\!=\!\!\overline{1,N}$. We notice that $\|\theta_i\|_\dn\!\!\to\!\!{0}$, $\forall i\!\!=\!\!\overline{1,N}$, where $\theta_i\!\!=\!\!(\mathcal{L}_i\!\otimes\!I_n)e$,  implies $\|(\mathcal{L}\!\otimes\!I_n)e\|\!\!\to\!\!0$. Thus we have $\diag\{\|\theta_i\|_\dn^\mu I_n\}_{i=1}^{N}\!(\mathcal{L}\!\otimes\! I_n)e\!\!\to\!\!\textbf{0}$ with $\|(\mathcal{L}\!\otimes\!I_n)e\|\!\!\to\!\!0$. Hence we
conclude the continuity of the function $e\!\!\to\!\!\nu\diag\{\|\theta_i\|_\dn^\mu I_n\}_{i=1}^{N}\!(\mathcal{L}\!\otimes\! I_n)e$ at $e\!\!=\!\!\textbf{0}$. 
\vspace{-0.2cm}

Based on the analysis above, for $q\!\!=\!\!\textbf{0}$, the error equation \eqref{eq:error_eq_fin} is continuous on $e\!\!\in\!\!\R^{Nn}$. 

\vspace{-0.4cm}
\section{}\label{app:finite}
\vspace{-0.3cm}
Using formula $\eqref{eq:procannorm}$, we derive
\vspace{-0.3cm}\begin{equation*}\begin{aligned}
    \tfrac{d \|e\|_{\tdn}}{d t}\!\!=\!\!\tfrac{\|e\|_{\tdn}e^{\top}\!\tdn^{\top}\!(-\ln \|e\|_{\tdn}) P\tdn(-\ln \|e\|_{\tdn})\dot e}{e^{\top}\!\tdn^{\top}\!(-\ln \|e\|_{\tdn})PG_{\tdn}\tdn(-\ln \|e\|_{\tdn})e},
	\end{aligned}\vspace{-0.3cm}\end{equation*}
 where
\vspace{-0.3cm}\begin{equation*}\begin{aligned}
\dot{e} \!\!=&(\tilde A\!+\!\diag\{g(|\omega_i|)\}_{i=1}^{N}\tilde H \tilde C\!-\!\nu\diag\{\|\theta_i\|_\dn^\mu I_n\}_{i=1}^{N}(\mathcal{L}\!\otimes\! I_n))e\\&\!+\!\Gamma(t,\hat x, x).    
\end{aligned}
\vspace{-0.3cm}\end{equation*}
Thus we have 
\vspace{-0.3cm}\begin{equation*}\begin{aligned}
&\|e\|_{\tdn}\tdn(-\ln \|e\|_{\tdn})\dot{e} \!\!=\!\!\|e\|_{\tdn}\tdn(-\ln \|e\|_{\tdn})\left(\tilde A\!+\!\diag\{g(|\omega_i|)\}_{i=1}^{N}\tilde H \tilde C\!\right.\\&\left.-\!\nu\diag\{\|\theta_i\|_\dn^\mu I_n\}_{i=1}^{N}(\mathcal{L}\!\otimes\! I_n)\right)e\!+\!\|e\|_{\tdn}\tdn(-\ln \|e\|_{\tdn})\Gamma(t,\hat x, x),
\end{aligned}\vspace{-0.3cm}
\end{equation*}
where 
$\|e\|_{\tdn}\tdn(-\ln \|e\|_{\tdn})\tilde A\!\!=\!\!\|e\|_{\tdn}^{\mu+1}\tilde A\tdn(-\ln \|e\|_{\tdn})$. 
Besides, 
\vspace{-0.3cm}\begin{equation*}
\begin{aligned}
\|e\|_{\tdn}\tdn(-\!\!\ln \|e\|_{\tdn})&\diag\{g(|\omega_i|)\}_{i=1}^{N}\tilde H \tilde C\\\!\!=&
\diag\{\|e\|_{\tdn}\dn(-\!\!\ln \|e\|_{\tdn})g(|\omega_i|) H_i  C_i\}_{i=1}^{N},\\
\end{aligned}
\vspace{-0.3cm}\end{equation*}
where
\vspace{-0.3cm}\begin{equation*}
\begin{aligned}
&\|e\|_{\tdn}\dn(-\!\!\ln \|e\|_{\tdn})g(|\omega_i|) H_i  C_i\\
\!\!=&\|e\|_{\tdn}\dn(-\!\!\ln \|e\|_{\tdn})\exp(\mu(G_0\!\!+\!\!I_n)\ln\!|C_ie_i|)H_i  C_i\\
\!\!=&\exp(\ln\|e\|_{\tdn})
\\&\times\!\!\exp(-\!\!\ln \|e\|_{\tdn}(\mu G_0\!\!+\!\!I_n))\exp(\mu(G_0\!\!+\!\!I_n)\ln\!|C_ie_i|)H_i  C_i\\
\!\!=&\exp(-\!\mu G_0 \ln \|e\|_{\tdn})
\\&\times\!\!\exp(\mu(G_0\!\!+\!\!I_n)\ln\!|C_i\dn(\ln \|e\|_{\tdn})\dn(-\!\!\ln \|e\|_{\tdn})e_i|)H_i  C_i\\
\!\!=&\exp(-\!\mu G_0 \ln \|e\|_{\tdn})
\\&\times\!\!\exp(\mu(G_0\!\!+\!\!I_n)\ln\!|\|e\|_{\tdn}C_i \dn(-\!\!\ln \|e\|_{\tdn})e_i|)H_i  C_i\\
\!\!=&\exp(-\!\mu G_0 \ln \|e\|_{\tdn})
\\&\times\!\!\exp(\mu(G_0\!\!+\!\!I_n)(\ln\!\|e\|_{\tdn}\!\!+\!\!\ln|C_i \dn(-\!\!\ln \|e\|_{\tdn})e_i|))H_i  C_i\\
\!\!=&\exp(-\!\mu G_0 \ln \|e\|_{\tdn})\exp(\mu(G_0\!\!+\!\!I_n)\ln\!\|e\|_{\tdn})\\&\!\!\times\!\exp(\mu(G_0\!\!+\!\!I_n)\ln|C_i \dn(-\!\!\ln \|e\|_{\tdn})e_i|)H_i  C_i\\
\!\!=&\|e\|_{\tdn}^{\mu}\exp(\mu(G_0\!\!+\!\!I_n)\ln|C_i \dn(-\!\!\ln \|e\|_{\tdn})e_i|)H_i  C_i\\
\!\!=&\|e\|_{\tdn}^{\mu+1}\exp(\mu(G_0\!\!+\!\!I_n)\ln|C_i \dn(-\!\!\ln \|e\|_{\tdn})e_i|)H_i  \|e\|_{\tdn}^{-1}C_i\\
\!\!=&\|e\|_{\tdn}^{\mu+1}\!\!\!\exp(\mu(G_0\!\!+\!\!I_n)\ln|C_i \dn(-\!\!\ln \|e\|_{\tdn})e_i|)H_i  C_i\dn(-\!\!\ln \|e\|_{\tdn}),
\end{aligned}
\vspace{-0.3cm}\end{equation*}
then 
\vspace{-0.3cm}\begin{equation*}
\begin{aligned}
&\|e\|_{\tdn}\tdn(-\!\!\ln \|e\|_{\tdn})\diag\{g(|\omega_i|)\}_{i=1}^{N}\tilde H \tilde C\\\!\!=&
\|e\|_{\tdn}^{\mu+1}
\\&\times\!\!\diag\{\exp(\mu(G_0\!\!+\!\!I_n)\!\!\ln\!\!|C_i \dn(-\!\!\ln \|e\|_{\tdn})e_i|)\}_{i=1}^{N}\!\tilde H \tilde C\tdn(-\!\!\ln \|e\|_{\tdn}).
\end{aligned}
\vspace{-0.3cm}\end{equation*}
Moreover, 
\vspace{-0.3cm}\begin{equation*}\begin{aligned}
&\nu\|e\|_{\tdn}\tdn(-\ln \|e\|_{\tdn})\diag\{\|\theta_i\|_\dn^\mu I_n\}_{i=1}^{N}(\mathcal{L}\!\otimes\! I_n)\\
=\!&\nu\|e\|_{\tdn}\diag\{\|(\mathcal{L}_i\!\otimes\! I_n)\tdn(\ln \|e\|_{\tdn})\tdn(-\ln \|e\|_{\tdn})e\|_\dn^\mu I_n\}_{i=1}^{N}\\&\times\!\tdn(-\ln \|e\|_{\tdn})(\mathcal{L}\!\otimes\! I_n).
\end{aligned}
\vspace{-0.3cm}\end{equation*}
Taking into account that $\tdn\!\!=\!\!I_N\!\otimes\!\dn$, and $(A\!\otimes\!B)(C\!\otimes\!D)\!\!=\!\!AC\!\otimes\!\!BD$, then we have $(\mathcal{L}_i\!\otimes\! I_n)\tdn(\ln \|e\|_{\tdn})\!\!=\!\!(\mathcal{L}_i\!\otimes\! I_n)(I_N\!\otimes\!\dn(\ln \|e\|_{\tdn}))\!\!=\!\!(\mathcal{L}_i\!\otimes\! \dn(\ln \|e\|_{\tdn}))\!\!=\!\!(1\!\otimes\!\dn(\ln \|e\|_{\tdn}))(\mathcal{L}_i\!\otimes\! I_n)\!\!=\!\!\dn(\ln \|e\|_{\tdn})(\mathcal{L}_i\!\otimes\! I_n)$ since $\mathcal{L}_i$ is the $i_{th}$ row of the Laplacian matrix. Meanwhile, $\tdn(-\ln \|e\|_{\tdn})(\mathcal{L}\!\otimes\! I_n)\!\!=\!\!(\mathcal{L}\!\otimes\! I_n)\tdn(-\ln \|e\|_{\tdn})$. Thus we have 
\vspace{-0.3cm}\begin{equation*}\begin{aligned}
&\nu\|e\|_{\tdn}\tdn(-\ln \|e\|_{\tdn})\diag\{\|\theta_i\|_\dn^\mu I_n\}_{i=1}^{N}(\mathcal{L}\!\otimes\! I_n)\\
=\!&\nu\|e\|_{\tdn}\diag\{\|\dn(\ln \|e\|_{\tdn})(\mathcal{L}_i\!\otimes\! I_n)\tdn(-\ln \|e\|_{\tdn})e\|_\dn^\mu I_n\}_{i=1}^{N}\\&\!\times\!(\mathcal{L}\!\otimes\! I_n)\tdn(-\ln \|e\|_{\tdn})\\
=\!&\nu\|e\|_{\tdn}^{1+\mu}\diag\{\|(\mathcal{L}_i\!\otimes\! I_n)\tdn(-\ln \|e\|_{\tdn})e\|_\dn^\mu I_n\}_{i=1}^{N}
\\&\times\!\!(\mathcal{L}\!\otimes\! I_n)\tdn(-\ln \|e\|_{\tdn}).
\end{aligned}
\vspace{-0.3cm}\end{equation*}
Based on the calculation above, we have  \eqref{eq:fin_error}.

\vspace{-0.4cm}
\section{}\label{app:muto0}
\vspace{-0.3cm}
Let the function $\sigma_1\!\!:\!\!\R\!\times\!\R^{p}\!\!\to\!\!\R^{Nn}$ defined as 
\vspace{-0.3cm}
\begin{equation*} \sigma_1(\mu,\tilde Cz)\!\!=\!\!(\tilde D(\tilde Cz)\!-\!I_{Nn})\tilde H\tilde Cz\!\!=\!\!\tilde D(\tilde Cz)\tilde H\tilde Cz\!\!-\!\!\tilde H\tilde Cz,
\vspace{-0.3cm}\end{equation*}
for $\tilde Cz\!\!\neq\!\!\textbf{0}$ and $\sigma_1(\mu,\textbf{0})\!\!=\!\!\textbf{0}$, where $\tilde D(\tilde Cz)\!\!=\!\!\diag\{g(|C_iz_i|)\}_{i=1}^{N}$, $g(|C_iz_i|)\!\!=\!\!\exp(\mu(G_0\!+\!I_n)\ln\!|C_iz_i|)$, $z\!\!=\!\!(z_1^{\top}\!\!,\dots,z_N^{\top})^{\top}\!\!\!=\!\tilde\dn(-\!\ln\!\|e\|_{\tilde\dn})e$.
It is clear that function $\sigma_1(\mu,\tilde Cz)$ is continuously differentiable on $(-1/{\tilde n},+\infty)\!\times\!(\R^p\backslash\{\textbf{0}\})$. Let us show it is continuously differentiable on $(-1/{\tilde n },+\infty)\!\times\!\R^p$ as well. 
On one hand, for any $\epsilon\!\!\in\!\!\R$ satisfying $\max(-\mu \tilde{n},0)\!\!\leq\!\! \epsilon \!\!<\!\!1$ we have
$\epsilon I_n\!+\!\mu(G_0\!+\!I_n)$ is anti-Hurwitz,
and taking into account 
\vspace{-0.3cm}\begin{equation*}\begin{aligned}
\ln\!|C_i&z_i|\exp(\mu(G_0\!+\!I_n)\ln\!|C_iz_i|)H_iC_iz_i\\=&\exp((\epsilon I_n\!\!+\!\!\mu(G_0\!\!+\!\!I_n))\ln\!|C_iz_i|)H_i\tfrac{C_iz_i\ln\!|C_iz_i|}{|C_iz_i|^\epsilon}, \; i\!\!=\!\!\overline{1,N}
\end{aligned}\vspace{-0.3cm}\end{equation*}
in which $\exp((\epsilon I_n\!\!+\!\!\mu(G_0\!\!+\!\!I_n))\ln\!|C_iz_i|)\!\!\to\!\!\textbf{0}$ as $|C_iz_i|\!\to\!0$, $\forall i\!\!=\!\!\overline{1,N}$ and $\tfrac{C_iz_i\ln\!|C_iz_i|}{|C_iz_i|^\epsilon}\!\!\to\!\!\textbf{0}$ as $|C_iz_i|\!\to\!0$, $\forall i\!\!=\!\!\overline{1,N}$.
Thus we say 
\vspace{-0.3cm}\begin{equation*}
\left(
\begin{smallmatrix}
...\\
\\
\ln\!|C_iz_i|\exp(\mu(G_0\!+\!I_n)\ln\!|C_iz_i|)H_iC_iz_i\\
\\
...\end{smallmatrix}\right)\!\!\to\!\!\textbf{0}\;as\;|C_iz_i|\!\!\to\!\!0,\;\forall i\!\!=\!\!\overline{1,N}.
\vspace{-0.3cm}\end{equation*}
On the other hand, $|\tilde Cz|\!\!\to\!\!0$ if and only if $|C_iz_i|\!\!\to\!\!0$, $\forall i\!\!=\!\!\overline{1,N}$.
Hence, we conclude 
\vspace{-0.3cm}\begin{equation*}
\tfrac{\partial\sigma_1(\mu,\tilde Cz)}{\partial\mu}\!\!=\!\!\left(
\begin{smallmatrix}
...\\
\\
(G_0\!+\!I_n)\ln\!|C_iz_i|\exp(\mu(G_0\!+\!I_n)\ln\!|C_iz_i|)H_iC_iz_i\\
\\
...\end{smallmatrix}\right)\!\!\to\!\!\textbf{0}\;as\;|\tilde Cz|\!\!\to\!\!0.
\vspace{-0.3cm}\end{equation*}
and
\vspace{-0.3cm}\begin{equation*}\begin{aligned}
&\sigma_1(\mu,\tilde Cz)\\
&\!\!=\!\!\left(
\begin{smallmatrix}
...\\
\\
\left(\exp((\epsilon I_n+\mu(G_0\!+\!I_n))\ln\!|C_iz_i|)\!-\!|C_iz_i|^\epsilon I_n\right)H_i\tfrac{C_iz_i}{|C_iz_i|^\epsilon}\\
\\
...\end{smallmatrix}\right)\!\!\!\to\!\!\textbf{0}\;as\;|\tilde Cz|\!\!\to\!\!0.
\end{aligned}\vspace{-0.3cm}\end{equation*}
The latter means 
$\sigma_1(\mu,\tilde Cz)$ is continuously differentiable on $(-1/{\tilde n},+\infty)\!\times\!\R^{p}$. Notice that for $z$ from the unit sphere one holds $|\tilde Cz|\!\!\in\!\![0,|\tilde CP^{-\frac{1}{2}}|]$, $P\!\!=\!\!I_N\!\!\otimes\!\!P_a$. Mean Value Theorem gives
\vspace{-0.3cm}\begin{equation*}\begin{aligned}
\|\sigma_1\!(\mu,\tilde Cz)\|_P^2\!\!=\!\!\|\sigma_1\!(0,\tilde Cz)\|_P^2\!+\!2|\mu|\sigma_1^\top\!\!(\tilde\mu,\tilde Cz)P\tfrac{\partial\sigma_1\!(\tilde \mu,\tilde Cz)}{\partial \tilde \mu}|_{\tilde\mu\!\in\![-|\mu|,|\mu|]}.
\end{aligned}\vspace{-0.3cm}\end{equation*}
Since $\|\sigma_1\!(0,\tilde Cz)\|_P^2\!\!=\!\!0$, then 
$\|\sigma_1(\mu,\tilde Cz)\|_P^2\!\!\leq\!\!2|\mu|\vartheta_1$, where 
\vspace{-0.3cm}$$
\vartheta_1\!\!:=\!\!\sup_{|\tilde\mu|\!
\leq\!|\mu|, |\tilde Cz|\!\leq\!|\tilde CP^{-{1}/{2}}|}|\sigma_1^\top\!\!(\tilde \mu,\tilde Cz)P\tfrac{\partial\sigma_1\!(\tilde \mu,\tilde Cz)}{\partial \tilde \mu}|\!\!<\!\!+\infty,
\vspace{-0.3cm}$$
and $\vartheta_1\!\!\to\!\!0$ with $\mu\!\!\to\!\!0$. Hence, $\sigma_1(\mu,\tilde C z)\!\!\to\!\!\textbf{0}$ as $\mu\!\!\to\!\!0$ uniformly on $z$ from the unit sphere.
\vspace{-0.2cm}

Let the function $\sigma_2\!\!:\!\!\R\!\times\!\R^{Nn}\!\!\to\!\!\R^{Nn}$ defined as
\vspace{-0.3cm}\begin{equation*}
  \begin{aligned}
      \sigma_2(\mu,(\mathcal{L}\!\otimes\! I_n)z)\!\!=\!\!(I_{Nn}\!-\!\tilde \Theta(z))(\mathcal{L}\!\otimes\! I_n)z,
  \end{aligned}
\vspace{-0.3cm}\end{equation*} 
where $\tilde \Theta(z)\!\!=\!\!\diag\{\|(\mathcal{L}_i\!\otimes\! I_n)z\|_\dn^\mu I_n\}_{i=1}^{N}$.
We show that this function is continuously differentiable on $(-\check \mu,+\check\mu)\!\times\!\R^{Nn}$, $0\!\!<\!\!\check\mu\!\!<\!\!1/{\tilde n}$ is small enough. 
On the one hand, 
for $i\!\!=\!\!\overline{1,N}$,
\vspace{-0.3cm}\begin{equation*}
    \begin{aligned}
    \|(\mathcal{L}_i\!\otimes\! I_n)z\|_\dn^{\mu}\ln(\|(\mathcal{L}_i&\otimes I_n)z\|_\dn)(\mathcal{L}_i\otimes I_n)z\!\!=\!\!\\
        &\|(\mathcal{L}_i\!\otimes\! I_n)z\|_\dn^{\mu+\varphi}\tfrac{\ln(\|(\mathcal{L}_i\otimes I_n)z\|_\dn)(\mathcal{L}_i\otimes I_n)z}{\|(\mathcal{L}_i\otimes I_n)z\|_\dn^\varphi},
    \end{aligned}
\vspace{-0.3cm}\end{equation*}
with $\varphi\!\!\in\!\!\R$ satisfying $1/\tilde n\!\!<\!\!\varphi\!\!<\!\!1$, then
we have $\forall \varphi\!\!\in\!\!(1/\tilde n,1)$,
\vspace{-0.3cm}\begin{equation*}
    \|(\mathcal{L}_i\!\otimes\! I_n)z\|_\dn^{\mu+\varphi}\ln(\|(\mathcal{L}_i\otimes I_n)z\|_\dn)\!\!\to\!\!{0}\;as\;\|(\mathcal{L}_i\!\otimes\! I_n)z\|\!\!\to\!\!0.
\vspace{-0.3cm}\end{equation*}
Besides, $\|(\mathcal{L}_i\!\otimes\! I_n)z\|\!\!\to\!\!0$ implies that $\|(\mathcal{L}_i\!\otimes\! I_n)z\|\!\!\in\!\!B(1)$, thus we have $\|(\mathcal{L}_i\!\otimes\! I_n)z\|_\dn^{\lambda_{\min}(G_\dn)}\!\!\geq\!\!\|(\mathcal{L}_i\!\otimes\! I_n)z\|_{P_a}$, therefore
\vspace{-0.3cm}\begin{equation*}
    \tfrac{\|(\mathcal{L}_i\otimes I_n)z\|_{P_a}}{\|(\mathcal{L}_i\otimes I_n)z\|_\dn^\varphi}\!\!\leq\!\!
    \tfrac{\|(\mathcal{L}_i\otimes I_n)z\|_{P_a}}{\|(\mathcal{L}_i\otimes I_n)z\|_{P_a}^{\frac{\varphi}{\lambda_{\min}(G_\dn)}}}\!\!=\!\!\|(\mathcal{L}_i\otimes I_n)z\|_{P_a}^{\frac{\lambda_{\min}(G_\dn)\!-\!\varphi}{\lambda_{\min}(G_\dn)}},
\vspace{-0.3cm}\end{equation*}
Since $\mu\!\!\in\!\!(-\check \mu,+\check\mu)$ is close enough to zero, then $G_\dn\!\!=\!\!\mu G_0\!\!+\!\!I_n\!\!\to\!\!I_    n$, and $\lambda_{\min}(G_\dn)\!\!\to\!\!1$, thus $\exists 1/\tilde n\!\!<\!\!\varphi\!\!<\!\!1$ such that $\varphi\!\!<\!\!\lambda_{\min}(G_\dn)$, then we have $\tfrac{(\mathcal{L}_i\otimes I_n)z}{\|(\mathcal{L}_i\otimes I_n)z\|_\dn^\varphi}\!\!\to\!\!\textbf{0}$ as $\|(\mathcal{L}_i\!\otimes\! I_n)z\|\!\!\to\!\!0$.
Thus we conclude 
$\|(\mathcal{L}_i\!\otimes\! I_n)z\|_\dn^{\mu}\ln(\|(\mathcal{L}_i\otimes I_n)z\|_\dn)(\mathcal{L}_i\otimes I_n)z\!\!\to\!\!\textbf{0}$,  as $\|(\mathcal{L}_i\!\!\otimes\!\! I_n)z\|\!\!\to\!\!{0}$, $\forall i\!\!=\!\!\overline{1,N}$.
On the other hand, $\|(\mathcal{L}\!\otimes\! I_n)z\|\!\!\to\!\!0$ is equivalent to $\|(\mathcal{L}_i\!\otimes\! I_n)z\|\!\!\to\!\!0$, $\forall i\!\!=\!\!\overline{1,N}$.
Thus 
\vspace{-0.3cm}\begin{equation*}\begin{aligned}
&\tfrac{\partial\sigma_2(\mu,(\mathcal{L}\!\otimes\! I_n)z)}{\partial\mu}\!\!=\!\!\left(\begin{smallmatrix}
    \dots\\
    \\
    \|(\mathcal{L}_i\!\otimes\! I_n)z\|_\dn^{\mu}\ln(\|(\mathcal{L}_i\otimes I_n)z\|_\dn)(\mathcal{L}_i\otimes I_n)z\\
    \\
    \dots
\end{smallmatrix}\right)\!\!\!\to\!\!\textbf{0},
\end{aligned}\vspace{-0.3cm}\end{equation*}
as $\|(\mathcal{L}\!\otimes\! I_n)z\|\!\!\to\!\!0$. In addition,
\vspace{-0.3cm}\begin{equation*}\begin{aligned}
&\sigma_2(\mu,(\mathcal{L}\!\otimes\! I_n)z)\!\!=\!\!\!\left(\begin{smallmatrix}
    \dots\\
    \\
    (\|(\mathcal{L}_i\otimes I_n)z\|_\dn^{\varphi}\!-\!\|(\mathcal{L}_i\otimes I_n)z\|_\dn^{\mu+\varphi})\frac{(\mathcal{L}_i\otimes I_n)z}{\|(\mathcal{L}_i\otimes I_n)z\|_\dn^{\varphi}}\\
    \\
    \dots
\end{smallmatrix}\right)\!\!\!\to\!\!\textbf{0},
\end{aligned}\vspace{-0.3cm}\end{equation*}
as $\|(\mathcal{L}\!\otimes\! I_n)z\|\!\!\to\!\!0$.
The latter gives
$\sigma_2(\mu,(\mathcal{L}\!\otimes\! I_n)z)$ is continuously differentiable on $(-\check \mu,+\check\mu)\!\times\!\R^{Nn}$.
Reusing Mean Value Theorem has
\vspace{-0.3cm}\begin{equation*}\begin{aligned}
\|\sigma_2\!(\mu,(\mathcal{L}\!\otimes\! &I_n)z)\|_P^2\!\!\leq\!\!\|\sigma_2\!(0,(\mathcal{L}\!\otimes\! I_n)z)\|_P^2\!\\
&+\!2|\mu|\sigma_2^\top\!\!(\tilde\mu,(\mathcal{L}\!\otimes\! I_n)z)P\tfrac{\partial\sigma_2\!(\tilde \mu,(\mathcal{L}\!\otimes\! I_n)z)}{\partial \tilde \mu}|_{\tilde\mu\!\in\![-|\mu|,|\mu|]}.
\end{aligned}\vspace{-0.3cm}\end{equation*}
Taking into account that $\|\sigma_2\!(0,(\mathcal{L}\!\otimes\! I_n)z)\|_P^2\!\!=\!\!0$ and $|(\mathcal{L}\!\otimes\! I_n)z|\!\!\leq\!\!|(\mathcal{L}\!\otimes\! I_n)P^{-\frac{1}{2}}|$, then
$\|\sigma_2\!(\mu,(\mathcal{L}\!\otimes\! I_n)z)\|_P^2\!\!\leq\!\!2|\mu|\vartheta_2$, with 
\vspace{-0.3cm}\begin{equation*}\begin{aligned}
\vartheta_2\!\!:=\!\!\!\!\sup_{|\tilde\mu|\!
\leq\!|\mu|, |(\mathcal{L}\!\otimes\! I_n)z|\!\leq\!|(\mathcal{L}\!\otimes\! I_n)P^{-{1}/{2}}|}\!\!\!\!|\sigma_2^\top\!\!(\tilde \mu,(\mathcal{L}\!\otimes\! I_n)z)P\tfrac{\partial\sigma_2\!(\tilde \mu,(\mathcal{L}\!\otimes\! I_n)z)}{\partial \tilde \mu}|,
\end{aligned}\vspace{-0.3cm}\end{equation*}
and $\vartheta_2\!\!\to\!\!0$ as $\mu\!\!\to\!\!0$. Hence, $\sigma_2(\mu,(\mathcal{L}\!\otimes\! I_n)z)\!\!\to\!\!\textbf{0}$ as $\mu\!\!\to\!\!0$ uniformly on $z$ from the unit sphere. 
\vspace{-0.4cm}
\section{}\label{app:fix}
\vspace{-0.3cm}
Recall equation \eqref{eq:fix_error} we have 
\vspace{-0.3cm}\begin{equation*}\begin{aligned}
\|e\|_{\tdn_k}&z^\top\!\!\! P\tdn_k(-\ln \|e\|_{\tdn_k})\dot{e}\!\!=\!\!\|e\|_{\tdn_k}\!z^\top\!\!\! P\tdn_k\!(\!-\!\ln \!\!\|e\|_{\tdn_k})\Gamma(t,\hat x, x) \\&\!\!+\!\!\|e\|_{\tdn_k}z^\top\!\!\! P\tdn_k(-\ln \|e\|_{\tdn_k})\left(\tilde A\!+\!\diag\{g(|\omega_i|)\}_{i=1}^{N}\tilde H \tilde C\!\right.
\\&\left.\!\!-\!\tfrac{\nu}{2}\diag\{\textstyle{\sum}_k\!\|\theta_i\|_{\dn_k}^{\mu_k} I_n\}_{i=1}^{N}\!\!(\mathcal{L}\!\!\otimes\!\! I_n)\!\!\!\right)\!\!e,
\end{aligned}
\vspace{-0.3cm}\end{equation*}
where $z\!\!=\!\!(z_1^{\top},\dots,z_N^{\top})^{\top}\!\!\!\!\!=\!\!\tilde\dn_k(-\ln\!\|e\|_{\tilde\dn_k})e$. 
In the latter equation, 
$\|e\|_{\tdn_k}z^\top\!\!\! P\tdn_k(-\ln \|e\|_{\tdn_k})\tilde Ae\!\!=\!\!\|e\|_{\tdn_k}^{\mu_k+1}z^\top\!\!\! P\tilde Az$. 
Besides, 
\vspace{-0.3cm}\begin{equation*}
\begin{aligned}
&\|e\|_{\tdn_k}z^\top\!\!\! P\tdn_k(-\!\!\ln \|e\|_{\tdn_k})\diag\{g(|\omega_i|)\}_{i=1}^{N}\tilde H \tilde Ce\\
\!\!=&\|e\|_{\tdn_k}z^\top\!\!\! P\tdn_k(-\!\!\ln \|e\|_{\tdn_k})\tfrac{1}{2}\diag\{\textstyle{\sum}_k\!\!\exp(\mu_k(G_0\!\!+\!\!I_n)\ln\!|C_ie_i|)\}\tilde H \tilde Ce\\
\!\!=&\tfrac{1}{2}\|e\|_{\tdn_k}z^\top\!\!\! P\tdn_k(-\!\!\ln \|e\|_{\tdn_k})\diag\{\exp(\mu_k(G_0\!\!+\!\!I_n)\ln\!|C_ie_i|)\}\tilde H \tilde Ce\\
\!\!+&\!\!\tfrac{1}{2}\|e\|_{\tdn_k}z^\top\!\!\! P\tdn_k(-\!\!\ln \|e\|_{\tdn_k})\diag\{\exp(\mu_{\tilde k}(G_0\!\!+\!\!I_n)\ln\!|C_ie_i|)\}\tilde H \tilde Ce,
\end{aligned}
\vspace{-0.3cm}\end{equation*}
with $\{k\}\!\!\cup\!\!\{\tilde k\}\!\!=\!\!\{0\}\!\!\cup\!\!\{\infty\}$, and
\vspace{-0.3cm}\begin{equation*}
\begin{aligned}
&\tfrac{1}{2}\|e\|_{\tdn_k}z^\top\!\!\! P\tdn_k(-\!\!\ln \|e\|_{\tdn_k})\diag\{\exp(\mu_k(G_0\!\!+\!\!I_n)\ln\!|C_ie_i|)\}\tilde H \tilde Ce\\
=&\tfrac{1}{2}\|e\|_{\tdn_k}^{\mu_k+1}z^\top\!\!\! P\diag\{\exp(\mu_k(G_0\!\!+\!\!I_n)\ln\!|C_iz_i|)\}\tilde H \tilde Cz,
\end{aligned}
\vspace{-0.3cm}\end{equation*}
and
\vspace{-0.3cm}\begin{equation*}
\begin{aligned}
&\tfrac{1}{2}\|e\|_{\tdn_k}z^\top\!\!\! P\tdn_k(-\!\!\ln \|e\|_{\tdn_k})\diag\{\exp(\mu_{\tilde k}(G_0\!\!+\!\!I_n)\ln\!|C_ie_i|)\}\tilde H \tilde Ce
\\=&\tfrac{1}{2}\|e\|^{\mu_k+1}_{\tdn_k}z^\top\!\!\! P\|e\|^{-\mu_k}_{\tdn_k}\tdn_k(-\!\!\ln \|e\|_{\tdn_k})
\\&\times\!\!\diag\{\exp(\mu_{\tilde k}(G_0\!\!+\!\!I_n)\ln\!|C_i\dn_k(\ln \|e\|_{\tdn_k})z_i|)\}\tilde H \tilde C\tdn_k(\ln \|e\|_{\tdn_k})z.
\end{aligned}
\vspace{-0.3cm}\end{equation*}
We notice that 
\vspace{-0.3cm}\begin{equation*}
    \begin{aligned}
&\lim_{\|e\|_{\tdn_k}\!\!\to\!k}\sup_{z^\top\!\!\! P z\!=\!\!1}\|\|e\|^{-\mu_k}_{\tdn_k}\tdn_k(-\!\ln\! \|e\|_{\tdn_k})
\\&\times\!\!\diag\{\exp(\mu_{\tilde k}(G_0\!\!+\!\!I_n)\ln\!|C_i\dn_k(\!\ln \|e\|_{\tdn_k})z_i|)\}\tilde H \tilde C\tdn_k(\!\ln \!\|e\|_{\tdn_k})z\|\!\!=\!\!0.
    \end{aligned}
\vspace{-0.3cm}\end{equation*}
Thus we say 
\vspace{-0.3cm}\begin{equation*}
    \begin{aligned}
&\lim_{\|e\|_{\tdn_k}\!\!\!\!\to\!k}\!\!\tfrac{1}{2}\|e\|_{\tdn_k}\!\!z^\top\!\!\! P\tdn_k(-\!\!\ln \|e\|_{\tdn_k})\diag\{\exp(\mu_{\tilde k}(G_0\!\!+\!\!I_n)\ln\!|C_ie_i|)\}\tilde H \tilde Ce\\&<\!\!\tfrac{\rho}{6}\|e\|_{\tdn_k}^{1+\mu_k}.
 \end{aligned}
\vspace{-0.3cm}\end{equation*}
Therefore,
\vspace{-0.3cm}\begin{equation*}
\begin{aligned}
&\|e\|_{\tdn_k}z^\top\!\!\! P\tdn_k(-\!\!\ln \|e\|_{\tdn_k})\diag\{g(|\omega_i|)\}_{i=1}^{N}\tilde H \tilde Ce\\
\!\!=&\tfrac{1}{2}\|e\|_{\tdn_k}z^\top\!\!\! P\tdn_k(-\!\!\ln \|e\|_{\tdn_k})\diag\{\exp(\mu_k(G_0\!\!+\!\!I_n)\ln\!|C_ie_i|)\}\tilde H \tilde Ce\\
\!\!&+\!\!\tfrac{1}{2}\|e\|_{\tdn_k}z^\top\!\!\! P\tdn_k(-\!\!\ln \|e\|_{\tdn_k})\diag\{\exp(\mu_{\tilde k}(G_0\!\!+\!\!I_n)\ln\!|C_ie_i|)\}\tilde H \tilde Ce,\\
\!\!<&\!\!\tfrac{1}{2}\|e\|_{\tdn_k}^{\mu_k+1}\!\!z^\top\!\!\! P\diag\{\exp(\mu_k(G_0\!\!+\!\!I_n)\ln\!|C_iz_i|)\}_{i=1}^N\tilde H \tilde  Cz\!\!+\!\!\tfrac{\rho}{6}\|e\|_{\tdn_k}^{1+\mu_k}
\end{aligned}
\vspace{-0.3cm}\end{equation*}
in the $k$-limit. 
Moreover, 
\vspace{-0.3cm}\begin{equation*}\begin{aligned}
&\tfrac{\nu}{2}\|e\|_{\tdn_k}z^\top\!\!\! P\tdn_k(-\ln \|e\|_{\tdn_k})\diag\{(\|\theta_i\|_{\dn_k}^{\mu_k}\!\!+\!\!\|\theta_i\|_{\dn_{\tilde k}}^{\mu_{\tilde k}})I_n\}_{i=1}^{N}(\mathcal{L}\!\otimes\! I_n)e\\
=\!&\tfrac{\nu}{2}\|e\|^{\mu_k+1}_{\tdn_k}z^\top\!\!\! P\diag\{\|(\mathcal{L}_i\!\otimes\! I_n)z\|_{\dn_k}^{\mu_k}I_n\}_{i=1}^{N}(\mathcal{L}\!\otimes\! I_n)z\\ \!\!+\!&\tfrac{\nu}{2}\|e\|_{\tdn_k}z^\top\!\!\! P\tdn_k(-\ln \|e\|_{\tdn_k})\diag\{\|(\mathcal{L}_i\!\otimes\! I_n)e\|_{\dn_{\tilde k}}^{\mu_{\tilde k}})I_n\}_{i=1}^{N}(\mathcal{L}\!\otimes\! I_n)e,
\end{aligned}
\vspace{-0.3cm}\end{equation*}
in which 
\vspace{-0.3cm}\begin{equation*}\begin{aligned}
&\tfrac{\nu}{2}z^\top\!\!\! P\|e\|_{\tdn_k}\tdn_k(-\!\!\ln \|e\|_{\tdn_k})\diag\{\|(\mathcal{L}_i\!\otimes\! I_n)e\|_{\dn_{\tilde k}}^{\mu_{\tilde k}}\}_{i=1}^N\!\!(\mathcal{L}\!\!\otimes\!\! I_n)e\\
&=\!\!\tfrac{\nu}{2}\|e\|_{\tdn_k}^{\mu_k+1}z^\top\!\!\! P\|e\|_{\tdn_k}^{-\mu_k}\tdn_k(-\ln \|e\|_{\tdn_k})\\&\times\diag\{\|(\mathcal{L}_i\!\otimes\! I_n)\tdn_k(\ln \|e\|_{\tdn_k})z\|_{\dn_{\tilde k}}^{\mu_{\tilde k}}\}_{i=1}^N(\mathcal{L}\!\otimes\! I_n)\tdn_k(\ln \|e\|_{\tdn_k})z.
\end{aligned}
\vspace{-0.3cm}\end{equation*}
We notice that 
\vspace{-0.3cm}\begin{equation*}\begin{aligned}
&\lim_{\|e\|_{\tdn_k}\!\!\to\!k}\sup_{z^\top\!\!\! Pz\!=\!1}\|\|e\|_{\tdn_k}^{-\mu_k}\tdn_k(-\ln \|e\|_{\tdn_k})\\&\times\diag\{\|(\mathcal{L}_i\!\otimes\! I_n)\tdn_k(\ln \|e\|_{\tdn_k})z\|_{\dn_{\tilde k}}^{\mu_{\tilde k}}\}(\mathcal{L}\!\otimes\! I_n)\tdn_k(\ln \|e\|_{\tdn_k})z\|\!\!=\!\!0.
\end{aligned}
\vspace{-0.3cm}\end{equation*}
Thus we say 
\vspace{-0.3cm}\begin{equation*}\begin{aligned}
&\lim_{\|e\|_{\tdn_k}\!\!\!\!\to\!k}\tfrac{\nu}{2}\|e\|_{\tdn_k}z^\top\!\!\! P\tdn_k(-\!\!\ln \|e\|_{\tdn_k})\diag\{\|(\mathcal{L}_i\!\otimes\! I_n)e\|_{\dn_{\tilde k}}^{\mu_{\tilde k}}\}_{i=1}^N\!\!(\mathcal{L}\!\!\otimes\!\! I_n)e\\\!\!<&\tfrac{\rho}{6}\|e\|_{\tdn_k}^{1+\mu_k}.
\end{aligned}
\vspace{-0.3cm}\end{equation*}
Therefore, 
\vspace{-0.3cm}\begin{equation*}\begin{aligned}
&\tfrac{\nu}{2}\|e\|_{\tdn_k}z^\top\!\!\! P\tdn_k(-\ln \|e\|_{\tdn_k})\diag\{(\|\theta_i\|_{\dn_k}^{\mu_k}\!\!+\!\!\|\theta_i\|_{\dn_{\tilde k}}^{\mu_{\tilde k}})I_n\}_{i=1}^{N}(\mathcal{L}\!\otimes\! I_n)e\\
<\!&\tfrac{\nu}{2}\|e\|^{\mu_k+1}_{\tdn_k}z^\top\!\!\! P\diag\{\|(\mathcal{L}_i\!\otimes\! I_n)z\|_{\dn_k}^{\mu_k}I_n\}_{i=1}^{N}(\mathcal{L}\!\otimes\! I_n)z \!\!+\!\!\tfrac{\rho}{6}\|e\|_{\tdn_k}^{1+\mu_k},
\end{aligned}
\vspace{-0.3cm}\end{equation*}
in the $k$-limit. 
Based on the calculation above we have  \eqref{eq:dercnk}.

\vspace{-0.4cm}
\section{}\label{app:rob_fin}
\vspace{-0.6cm}
\vspace{-0.4cm}\begin{equation*}
\begin{aligned}
&\|e\|_{\tdn}\tdn(-\!\!\ln \|e\|_{\tdn})\diag\{g(|C_ie_i\!\!-\!\!q_{y,i}|)\}_{i=1}^{N}\tilde H (\tilde Ce\!\!-\!\!q_y)\\\!\!=&
\diag\{\|e\|_{\tdn}\dn(-\!\!\ln \|e\|_{\tdn})g(|C_ie_i\!\!-\!\!q_{y,i}|) H_i(  C_ie_i\!\!-\!\!q_{y_i})\}_{i=1}^{N},\\
\end{aligned}
\vspace{-0.3cm}\end{equation*}
where
\vspace{-0.3cm}\begin{equation*}
\begin{aligned}
&\|e\|_{\tdn}\dn(-\!\!\ln \|e\|_{\tdn})g(|C_ie_i\!\!-\!\!q_{y,i}|) H_i(  C_ie_i\!\!-\!\!q_{y_i})\\
=&\exp(-\mu G_0\ln\|e\|_{\tdn})
\\&\times\!\!g(|C_i\dn(\ln \|e\|_{\tdn})\dn(-\!\!\ln \|e\|_{\tdn})e_i\!\!-\!\!q_{y,i}|) H_i(  C_ie_i\!\!-\!\!q_{y_i})\\
=&\exp(-\mu G_0\ln\|e\|_{\tdn})g(\|e\|_{\tdn}|\epsilon_i|) H_i(  C_ie_i\!\!-\!\!q_{y_i})\\
=&\exp(-\mu G_0\ln\|e\|_{\tdn})
\exp(\mu(G_0\!+\!I_n)\ln(\|e\|_{\tdn}|\epsilon_i|)) H_i(  C_ie_i\!\!-\!\!q_{y_i})\\
=&\exp(-\mu G_0\ln\|e\|_{\tdn})\exp(\mu(G_0\!+\!I_n)\ln\|e\|_{\tdn})
\\&\times\!\!\exp(\mu(G_0\!+\!I_n) \ln(|\epsilon_i|)) H_i(  C_ie_i\!\!-\!\!q_{y_i})\\
=&\|e\|_{\tdn}^{\mu}\exp(\mu(G_0\!+\!I_n) \ln(|\epsilon_i|)) H_i(  C_ie_i\!\!-\!\!q_{y_i})\\
=&\|e\|_{\tdn}^{\mu+1}\exp(\mu(G_0\!+\!I_n) \ln(|\epsilon_i|)) H_i \epsilon_i,
\end{aligned}
\vspace{-0.3cm}\end{equation*}
with $\epsilon_i\!\!=\!\!C_iz_i\!\!-\!\!\|e\|_{\tdn}^{-1}q_{y,i}$.

\vspace{-0.4cm}
\section{}\label{app:rob_fin_second}
\vspace{-0.3cm}
Notice $|\epsilon|\!\!=\!\!|\tilde Cz\!-\!\|e\|_{\tilde\dn}^{-1}\!q_y|\!\!\leq\!\!|\tilde CP^{-1/2}|\!\!+\!\!\|e\|_{\tdn}^{-1}|q_{y}|$, and 
 $\exists \pi\!\!>\!\!{|\tilde CP^{-{1}/{2}}|}$, such that 
$|\epsilon|\!\!\leq\!\!|\tilde CP^{-1/2}|\!\!+\!\!\|e\|_{\tilde\dn}^{-1}\!|q_{y}|\!\!<\!\!\pi$  provided $\|e\|_{\tilde\dn}\!\!>\!\!\tfrac{\sqrt{p}}{\pi\!\!-\!\!|\tilde CP^{-{1}/{2}}|}\|q_y\|_{L_\infty}$. 
Since $|\epsilon|\!\!=\!\!\sqrt{\sum_{i=1}^p|\epsilon_i|^2}$, thus $|\epsilon_i|\!\!<\!\!|\epsilon|\!\!<\!\!\pi$, $\forall i\!\!=\!\!\overline{1,p}$.
Therefore, similar to Appendix \ref{app:muto0}, we have  $\sup_{|\epsilon_i|\!<\!\pi}(\exp(\mu(G_0\!\!+\!\!I_n)|\epsilon_i|)\!\!-\!\!I_n)H_i\epsilon_i\!\!\to\!\!\textbf{0}$ as $\mu\!\!\to\!\!0$, $\forall i\!\!=\!\!\overline{1,p}$. The latter implies $(\exp(\mu(G_0\!\!+\!\!I_n)|\epsilon_i|)\!\!-\!\!I_n)H_i\epsilon_i\!\!\to\!\!\textbf{0}$ as $\mu\!\!\to\!\!0$, $\forall i\!\!=\!\!\overline{1,p}$ uniformly.  Therefore, with $\mu$ sufficiently close to zero, 
$(\tilde D(\tilde Cz,\|e\|_{\tilde\dn}^{-1}\!q_y)\!\!-\!\!I_{Nn})\tilde H\epsilon\!\!\to\!\!\textbf{0}$ uniformly. So we say $z^\top\!\!\! P(\tilde D(\tilde Cz,\|e\|_{\tilde\dn}^{-1}\!q_y)\!\!-\!\!I_{Nn})\tilde H\epsilon\!\!<\!\!\frac{\rho}{9}$.

\vspace{-0.4cm}
\section{}\label{app:rob_fin_third}
\vspace{-0.3cm}
Since $\tilde\dn_{q}$ is generated by $G_{\tilde\dn_{q}}\!\!\!\!\!=\!\!I_N\!\!\otimes\!\!(\mu(G_0\!\!+\!\!I_n)\!\!+\!\!I_n)$, we have
\vspace{-0.3cm}\begin{equation*}
\|e\|_{\tilde\dn}^{-\mu}\!\!z^\top \!\!\!P\tilde\dn(-\!\ln\!\|e\|_{\tilde\dn})\tilde q_x\!\!
=\!\!z^\top \!\!\!P\tilde\dn_{q}(-\!\ln\!\|e\|_{\tilde\dn})\tilde q_x,
\vspace{-0.3cm}\end{equation*}
$G_{\tilde\dn_{q}}$ is anti-Hurwitz since $G_0$ is obtained from \eqref{eq:G_0} and $\mu\!\!>\!\!-1/{\tilde n}$.
Then, 
\vspace{-0.3cm}\begin{equation*}
    \begin{aligned}
z^\top \!\!\!P\tilde\dn_{q}(-\!\ln\!\|e\|_{\tilde\dn})\tilde q_x
&\!\!=\!\!z^\top \!\!\!P\tilde\dn_{q}(-\!\ln\!\tfrac{\|e\|_{\tilde\dn}}{\|\tilde q_x\|_{\tilde\dn_{q}}})\tilde\dn_{q}(-\!\ln\!\|\tilde q_x\|_{\tilde\dn_{q}})\tilde q_x\\
&\!\!\leq\!\!|P^{\frac{1}{2}}||\tilde\dn_{q}(-\ln\!\tfrac{\|e\|_{\tilde\dn}}{\|\tilde q_x\|_{\tilde\dn_{q}}})||P_{q}^{-\frac{1}{2}}|,
\end{aligned}\vspace{-0.3cm}\end{equation*}
 the canonical homogeneous norm $\|\cdot\|_{\tilde\dn_{q}}$ is induced by weighted Euclidean norm $\|\cdot\|_{P_{q}}$.
Then
\vspace{-0.4cm}\begin{equation*}
    \begin{aligned}
z^\top \!\!\!P\tilde\dn_{q}(-\!\ln\!\|e\|_{\tilde\dn})\tilde q_x
\!\!\leq\!\!|P^{\frac{1}{2}}||P_{q}^{-\frac{1}{2}}|\tfrac{\|\tilde q_x\|^{\lambda_m}_{\tilde\dn_{q}}}{\|e\|^{\lambda_m}_{\tilde\dn}},
\end{aligned}\vspace{-0.3cm}\end{equation*}
provided $\|e\|_{\tilde\dn}\!\!\geq\!\!\|\tilde q_x\|_{\tilde\dn_{q}}$,
where $\lambda_m\!\!=\!\!\lambda_{\min}(G_{\tilde\dn_{q}})\!\!>\!\!0$.
Thus we have  $z^\top \!\!\!P\tilde\dn_{q}(-\!\ln\!\|e\|_{\tilde\dn})\tilde q_x\!\!<\!\!\frac{\rho}{9}$ provided $\|e\|_{\tilde\dn}\!\!>\!\!\Upsilon_M\|\tilde q_x\|_{\tilde\dn_{q}}$, $\Upsilon_M\!\!=\!\!\max\{1,\xi_M^{-1}\}$.

\vspace{-0.4cm}
\section{}\label{app:rob_fix}
\vspace{-0.3cm}
Using formula $\eqref{eq:procannorm}$, we derive
\vspace{-0.3cm}\begin{equation*}\begin{aligned}
    \tfrac{d \|e\|_{\tdn_\infty}}{d t}\!\!=\!\!\tfrac{\|e\|_{\tdn_\infty}e^{\top}\!\tdn_\infty^{\top}\!(-\ln \|e\|_{\tdn_\infty}) P\tdn_\infty(-\ln \|e\|_{\tdn_\infty})\dot e}{e^{\top}\!\tdn_\infty^{\top}\!(-\ln \|e\|_{\tdn_\infty})PG_{\tdn_\infty}\tdn_\infty(-\ln \|e\|_{\tdn_\infty})e},
	\end{aligned}\vspace{-0.3cm}\end{equation*}
 where
\vspace{-0.3cm}\begin{equation*}\begin{aligned}
\dot{e} \!\!=\tilde Ae\!+\!\diag\{g(&|\omega_i|)\}_{i=1}^{N}\tilde H (\tilde Ce\!-\!q_y)\!+\!\Gamma(t,\hat x, x)\!-\!\tilde q_x\\&\!-\!\tfrac{\nu}{2}\diag\{(\|\theta_i\|_{\dn_0}^{\mu_0}\!\!+\!\!\|\theta_i\|_{\dn_\infty}^{\mu_\infty}) I_n\}_{i=1}^{N}(\mathcal{L}\!\otimes\! I_n)e,    
\end{aligned}
\vspace{-0.3cm}\end{equation*}
with $\omega_i\!\!=\!\!C_ie_i\!\!-\!\!q_{y,i}$.
Thus we have 
\vspace{-0.3cm}\begin{equation*}\begin{aligned}
&\|e\|_{\tdn_\infty}z^\top\!\!\! P\tdn_\infty(-\ln \|e\|_{\tdn_\infty})\dot{e} \!\!=\!\!\|e\|_{\tdn_\infty}z^\top\!\!\! P\tdn_\infty(-\ln \|e\|_{\tdn_\infty})\tilde Ae\\
&\!+\!\|e\|_{\tdn_\infty}z^\top\!\!\! P\tdn_\infty(-\ln \|e\|_{\tdn_\infty})\diag\{g(|\omega_i|)\}_{i=1}^{N}\tilde H (\tilde Ce\!-\!q_y)\\
&\!+\!\|e\|_{\tdn_\infty}z^\top\!\!\! P\tdn_\infty(-\ln \|e\|_{\tdn_\infty})(\Gamma(t,\hat x, x)\!-\!\tilde q_x)\\
&\!-\!\!\tfrac{\nu}{2}\|e\|_{\tdn_\infty}\!\!\!z^\top\!\!\! P\tdn_\infty\!(-\ln \!\!\|e\|_{\tdn_\infty})\diag\{(\|\theta_i\|_{\dn_0}^{\mu_0}\!\!+\!\!\|\theta_i\|_{\dn\infty}^{\mu_\infty}) I_n\}_{i=1}^{N}\!\!(\mathcal{L}\!\!\otimes\!\! I_n)e,
\end{aligned}
\vspace{-0.3cm}\end{equation*}
where $z\!\!=\!\!(z_1^{\top},\dots,z_N^{\top})^{\top}\!\!\!\!\!=\!\!\tilde\dn_\infty(-\ln\!\|e\|_{\tilde\dn_\infty})e$, $\mu_\infty\!\!>\!\!0$ and $\mu_0\!\!<\!\!0$. 
In the latter equation, 
we detail the calculation of the second term of the right-hand side, while the others are similar to the proof of Theorem \ref{thm:bilimit_stab} with $k\!\!=\!\!\infty$. 
\vspace{-0.3cm}\begin{equation*}
\begin{aligned}
&\|e\|_{\tdn_\infty}z^\top\!\!\! P\tdn_\infty(-\!\!\ln \|e\|_{\tdn_\infty})\diag\{g(|C_ie_i\!\!-\!\!q_{y,i}|)\}_{i=1}^{N}\tilde H (\tilde Ce\!\!-\!\!q_y)\\
\!\!=&\tfrac{1}{2}\|e\|_{\tdn_\infty}\!\!z^\top\!\!\! P\tdn_\infty(\!\!-\!\!\ln\!\! \|e\|_{\tdn_\infty})\!\!\\
&\times \diag\{\exp(\mu_\infty(G_0\!\!+\!\!I_n)\!\!\ln\!|C_ie_i\!\!-\!\!q_{y,i}|)\}\!\tilde H  (\tilde Ce\!\!-\!\!q_{y})\\
\!\!&+\!\!\tfrac{1}{2}\|e\|_{\tdn_\infty}\!\!z^\top\!\!\! P\tdn_\infty(\!\!-\!\!\ln \|e\|_{\tdn_\infty})\!\!\\
&\times\diag\{\exp(\mu_{0}(G_0\!\!+\!\!I_n)\!\!\ln\!|C_ie_i\!\!-\!\!q_{y,i}|)\}\!\tilde H(\tilde Ce\!\!-\!\!q_{y}),
\end{aligned}
\vspace{-0.3cm}\end{equation*}
and
\vspace{-0.3cm}\begin{equation*}
\begin{aligned}
&\tfrac{1}{2}\|e\|_{\tdn_\infty}\!\!z^\top\!\!\! P\tdn_\infty(\!\!-\!\!\ln\!\! \|e\|_{\tdn_\infty})\!\!\\
&\times \diag\{\exp(\mu_\infty(G_0\!\!+\!\!I_n)\!\!\ln\!|C_ie_i\!\!-\!\!q_{y,i}|)\}\!\tilde H  (\tilde Ce\!\!-\!\!q_{y})\\
=&\tfrac{1}{2}\|e\|_{\tdn_\infty}^{\mu_\infty\!+\!1}\!\!z^\top\!\!\! P\!\!\\
&\times \diag\{\exp(\mu_\infty(G_0\!\!+\!\!I_n)\!\!\ln\!|C_iz_i\!\!-\!\!\|e\|_{\tdn_\infty}^{-1}q_{y,i}|)\}\!\tilde H  (\tilde Cz\!\!-\!\!\|e\|_{\tdn_\infty}^{-1}q_{y}),
\end{aligned}
\vspace{-0.3cm}\end{equation*}
and
\vspace{-0.3cm}\begin{equation*}
\begin{aligned}
&\!\!\tfrac{1}{2}\|e\|_{\tdn_\infty}\!\!z^\top\!\!\! P\tdn_\infty(\!\!-\!\!\ln \|e\|_{\tdn_\infty})\!\!\\
&\times\diag\{\exp(\mu_{0}(G_0\!\!+\!\!I_n)\!\!\ln\!|C_ie_i\!\!-\!\!q_{y,i}|)\}\!\tilde H(\tilde Ce\!\!-\!\!q_{y})\\
=&\tfrac{1}{2}\|e\|_{\tdn_\infty}^{\mu_\infty+1}\!\!z^\top\!\!\! P\|e\|_{\tdn_\infty}^{-\mu_\infty}\tdn_\infty(\!\!-\!\!\ln \|e\|_{\tdn_\infty})\!\!\\
&\times\diag\{\exp(\mu_{0}(G_0\!\!+\!\!I_n)\!\!\ln\!|C_i\dn_\infty(\ln \|e\|_{\tdn_\infty})z_i\!\!-\!\!q_{y,i}|)\}\!\\&\times\tilde H(\tilde C\tdn_\infty(\ln \|e\|_{\tdn_\infty})z\!\!-\!\!q_{y})\\
\end{aligned}
\vspace{-0.3cm}\end{equation*}
We notice that for $q_{y}\!\!\in\!\!L^\infty(\R,\R^p)$,
\vspace{-0.3cm}\begin{equation*}
    \begin{aligned}
\lim_{\|e\|_{\tdn_\infty}\!\!\to\!\infty}\sup_{z^\top\!\!\! P z\!=\!\!1}&
\|\|e\|_{\tdn_\infty}^{-\mu_\infty}\tdn_\infty(\!\!-\!\!\ln \|e\|_{\tdn_\infty})\!\!\\
&\times\!\!\diag\{\exp(\mu_{0}(G_0\!\!+\!\!I_n)\!\!\ln\!|C_i\dn_\infty(\ln \|e\|_{\tdn_\infty})z_i\!\!-\!\!q_{y,i}|)\}\!\\
&\times\!\!\tilde H(\tilde C\tdn_\infty(\ln \|e\|_{\tdn_\infty})z\!\!-\!\!q_{y})\|\!\!=\!\!0.
    \end{aligned}
\vspace{-0.3cm}\end{equation*}
Thus we say 
\vspace{-0.3cm}\begin{equation*}
    \begin{aligned}
&\lim_{\|e\|_{\tdn_\infty}\!\!\to\!\infty}\!\!\tfrac{1}{2}\|e\|_{\tdn_\infty}\!\!z^\top\!\!\! P\tdn_\infty(\!\!-\!\!\ln \|e\|_{\tdn_\infty})\!\!\\
&\times\diag\{\exp(\mu_{0}(G_0\!\!+\!\!I_n)\!\!\ln\!|C_ie_i\!\!-\!\!q_{y,i}|)\}\!\tilde H(\tilde Ce\!\!-\!\!q_{y})\!\!<\!\!\tfrac{\rho}{9}\|e\|_{\tdn_\infty}^{1+\mu_\infty}.
 \end{aligned}
\vspace{-0.3cm}\end{equation*}
Then,
\vspace{-0.3cm}\begin{equation*}
\begin{aligned}
&\|e\|_{\tdn_\infty}z^\top\!\!\! P\tdn_\infty(-\!\!\ln \|e\|_{\tdn_\infty})\diag\{g(|\omega_i|)\}_{i=1}^{N}\tilde H (\tilde Ce\!\!-\!\!q_y)\\
\!\!=&\tfrac{1}{2}\|e\|_{\tdn_\infty}\!\!z^\top\!\!\! P\tdn_\infty(\!\!-\!\!\ln\!\! \|e\|_{\tdn_\infty})\! \diag\{\exp(\mu_\infty(G_0\!\!+\!\!I_n)\!\!\ln\!|\omega_i|)\}\!\tilde H  (\tilde Ce\!\!-\!\!q_{y})\\
\!\!+&\!\!\tfrac{1}{2}\|e\|_{\tdn_\infty}\!\!z^\top\!\!\! P\tdn_\infty(\!\!-\!\!\ln \|e\|_{\tdn_\infty})\!\diag\{\exp(\mu_{0}(G_0\!\!+\!\!I_n)\!\!\ln\!|\omega_i|)\}\!\tilde H(\tilde Ce\!\!-\!\!q_{y}),\\
\!\!<&\tfrac{1}{2}\|e\|_{\tdn_\infty}^{\mu_\infty+1}\!\!z^\top\!\!\! P\tilde D_\infty(\tilde Cz,\|e\|_{\tilde\dn_\infty}^{-1}\!q_y)\tilde H (\tilde  Cz\!\!-\!\!\|e\|_{\tdn_\infty}^{-1}q_{y})\!\!+\!\!\tfrac{\rho}{9}\|e\|_{\tdn_\infty}^{1+\mu_\infty},
\end{aligned}
\vspace{-0.3cm}\end{equation*}
with $\tilde D_\infty(\tilde Cz,\|e\|_{\tilde\dn_\infty}^{-1}\!q_y)\!\!=\!\!\diag\{\exp(\mu_\infty(G_0\!\!+\!\!I_n)\ln\!|C_iz_i\!\!-\!\!\|e\|_{\tdn_\infty}^{-1}q_{y,i}|)\}_{i=1}^N$. 
\vspace{-0.2cm}

Following Theorem \ref{thm:bilimit_stab}, we estimate 
\vspace{-0.3cm}\begin{equation*}
\begin{aligned}
&\tfrac{\nu}{2}\|e\|_{\tdn_\infty}\!\!\!z^\top\!\!\! P\tdn_\infty\!(-\ln \!\!\|e\|_{\tdn_\infty})\diag\{(\|\theta_i\|_{\dn_0}^{\mu_0}\!\!+\!\!\|\theta_i\|_{\dn\infty}^{\mu_\infty}) I_n\}_{i=1}^{N}\!\!(\mathcal{L}\!\!\otimes\!\! I_n)e
\\
<\!&\tfrac{\nu}{2}\|e\|^{\mu_\infty+1}_{\tdn_\infty}\!\!z^\top\!\!\! P\diag\{\|(\mathcal{L}_i\!\otimes\! I_n)z\|_{\dn_\infty}^{\mu_\infty}I_n\}_{i=1}^{N}(\mathcal{L}\!\otimes\! I_n)z \!\!+\!\!\tfrac{\rho}{6}\|e\|_{\tdn_\infty}^{1+\mu_\infty}.
\end{aligned}
\vspace{-0.3cm}\end{equation*}
Based on the calculation above we have \eqref{eq:rob_fix_estimation}.

\bibliographystyle{agsm}
\bibliography{root}
\vspace{-0.5cm}

\end{document}